\documentclass[a4paper,12pt]{article}
\usepackage[T1]{fontenc}
\usepackage[utf8]{inputenc}
\usepackage[english]{babel}
\newcommand\Dcyr{{\mathfrak D}}
\newcommand\dcyr{{\mathfrak d}}
\newcommand\Lcyr{{\mathfrak L}}
\usepackage[margin=2.5cm]{geometry}
\usepackage{amsmath,amsthm,amsfonts,stmaryrd,amssymb}
\usepackage{remyrefs}
\usepackage[all]{xy}

\DeclareMathAlphabet{\mathpzc}{T1}{pzc}{m}{it}
\def\leq{\leqslant}
\def\geq{\geqslant}
\def\eps{\varepsilon}
\def\phi{\varphi}
\renewcommand\emptyset\varnothing

\newcommand\Z{\mathbb{Z}}

\newcommand\C{\mathbb{C}}

\newcommand\pa[1]{\ensuremath{\left(#1\right)}}
\newcommand\set[1]{\ensuremath{\left\{#1\right\}}}
\newcommand\abs[1]{\ensuremath{\left\lvert#1\right\rvert}}

\newcommand\irange[1]{\ensuremath{\left\llbracket#1\right\rrbracket}}

\newcommand\pair[1]{\ensuremath{\left\langle#1\right\rangle}}

\newcommand\vect{\boldsymbol}

\DeclareMathOperator{\rank}{rk}
\DeclareMathOperator{\Tr}{Tr}
\DeclareMathOperator{\Hom}{Hom}
\DeclareMathOperator{\Ext}{Ext}
\DeclareMathOperator{\Tor}{Tor}

\renewcommand{\O}{\mathcal O}

\DeclareMathOperator{\Pic}{Pic}
\DeclareMathOperator{\Jac}{Jac}
\DeclareMathOperator{\Mor}{Mor}

\DeclareMathOperator{\Homrond}{\mathpzc{Hom}}

\DeclareMathOperator{\Quot}{Quot}

\DeclareMathAlphabet{\mathss}{OT1}{cmss}{b}{n}

\DeclareMathOperator{\gr}{gr}
\DeclareMathOperator{\Gr}{Gr}
\DeclareMathOperator{\Vdm}{Vdm}

\def\arXiv#1{%
    \href{http://arxiv.org/abs/#1}{arXiv:#1}%
}

\theoremstyle{plain}
\newtheorem{thm}{Theorem}[section]

\newtheorem{defn}[thm]{Definition}

\newtheorem{prop}[thm]{Proposition}
\newtheorem{lemma}[thm]{Lemma}
\newtheorem{cor}[thm]{Corollary}

\theoremstyle{remark}
\newtheorem*{rem}{Remark}

\title{Rank-level duality for \\ conformal blocks of $SL_r$ and $GL_r$}
\author{Rémy Oudompheng}
\date\today

\begin{document}

\maketitle

\begin{abstract}
  This work aims at generalising the rank-level duality (also known as
  \emph{strange duality}) proved by Belkale \cite{Belkale} for generic
  curves and Marian and Oprea \cite{MarianOprea} for every smooth curve to
  the case of spaces of conformal blocks related to moduli spaces of
  parabolic bundles on a smooth projective curve.
\end{abstract}

\tableofcontents

\newpage

\section*{Introduction}

\def\Xp{{X,\vect p}}
\def\Upar#1{\mathcal U_{X,\vect p}(#1, \vect\lambda)}
\def\SUpar{\mathcal S\Upar}
\def\Lpar#1{\mathcal L^{#1, \vect\lambda}}
\def\Utpar#1{\mathcal U_{X,\vect p}(#1, \vect{\lambda^T})}
\def\SUtpar{\mathcal S\Utpar}
\def\Ltpar#1{\mathcal L^{#1, \vect{\lambda^T}}}
\def\Uppar#1{\mathcal U_{X,\vect p}(#1, \vect{\abs{\lambda}})}
\def\SUppar{\mathcal S\Uppar}
\def\Lppar#1{\mathcal L^{#1, \vect{\abs{\lambda}}}}

Conformal blocks are vector spaces which appear for example in the
Conformal Field Theory of WZW models \cite{Witten-QFT}: mathematically
speaking, they can be defined as sections of determinant line bundles (and
their powers) on moduli spaces of vector bundles. The strange duality
conjecture, which appears as a mathematical statement in \cite{DonagiTu},
asserts that these spaces are dual to each other when switching rank 
(of vector bundles) and level (the power of the line bundle). 

The dimension of these spaces can be computed by the Verlinde formula
\citelist{\cite{Verlinde} \cite{Beauville-verlinde}} and can be checked to
be equal. The strange duality conjecture was proved recently by Belkale
\citelist{\cite{Belkale} \cite{Belkale-parallel}} for rational nodal curves
(and thus for generic curves) and extended to arbitrary smooth curves of
genus $g \geq 2$ by Marian and Oprea \cite{MarianOprea}. The proofs are
based on the isomorphism between the Verlinde algebra and the quantum
cohomology of the Grassmannian, which is explained by a paper by Witten
\cite{Witten} using quantum field theory language and arguments: this
translates into spectacular enumerative properties which were first used by
Belkale, then by Marian and Oprea in a slightly different fashion, to prove
the conjecture.

This text aims to explain how Marian and Oprea's proof can be generalised
in a straightforward way to parabolic conformal blocks on moduli spaces of
vector bundles with parabolic structures at marked points. Its goal is also
to give a geometric proof of the statements of Nakanishi and Tsuchiya
\cite{NakanishiTsuchiya}.

I would like to thank Arnaud Beauville for his remarks and corrections. 

\setcounter{section}{-1}
\section{Notations}

Throughout this text, $(\Xp)$ will denote a smooth projective curve over
$\C$ of genus $g$ with $n$ marked points, and $r$ and $l$ are
fixed positive integers. We denote by $D$ the finite set of marked points.

We will attach to the marked points $n$ Young diagrams denoted by
$\vect\lambda = (\lambda_p)_{p \in D}$: each of this diagram is a
nonincreasing sequence $\lambda = (a_1 \geq a_2 \geq \cdots \geq a_r)$ of
$r$ integers such that $0 \leq a_i \leq l$. These Young diagrams will often
represent dominant weights of the Lie algebra $\mathfrak{sl}_r$ (which can
be canonically associated to nonincreasing sequences of $r$ integers up to
translation) of level $\leq l$ (the level of a weight being the
difference between the first and the last number of the sequence). A Young
diagram can also be represented by a set of boxes in a $l \times r$
rectangle divided in $r$ rows whose lengths are the $a_i$'s. The set of
such Young diagrams is denoted by $YD_l$. 

The transpose of a diagram (represented by its Young diagram of height $r$
and width $l$) is written $\lambda^T$: it is the image of $\lambda$ under
the reflection which maps the horizontal and vertical axes to each other.
It can represent a dominant weight of $\mathfrak{sl}_l$ of level at most
$r$.

The conjugate diagram of $\lambda$ is defined by the sequence $\lambda^\ast
= (l-a_r, \ldots, l-a_1)$: its graphical representation consists of the
boxes in the complement of the diagram of $\lambda$ in the $l \times r$
rectangle. 

The \emph{size} $\abs\lambda$ of $\lambda$, is the number of boxes of its
Young diagram, which is $\sum a_i$ in the previous notations. If
$\vect\lambda$ if a family of weights, $\abs{\vect\lambda}$ is the sum of
their sizes.

If $\lambda = (a_1 \geq a_2 \geq \cdots \geq a_r)$, we will use the
notation $\lambda^{(i)} := a_i$ to refer to the components of $\lambda$ and 
the notation $Q_\lambda(X_1, \ldots, X_r)$ for the antisymmetric polynomial 
$\det (X_i^{a_j+r-j})_{i,j}$, $\Vdm := Q_0$ for the Vandermonde determinant 
and $S_\lambda(X_1, \ldots, X_r)$ for the Schur polynomials
$Q_\lambda/\Vdm$.

\section{Description of the strange duality}

\subsection{Parabolic structures on vector spaces}

Let $\mathfrak g := \mathfrak{sl}_r$ be the Lie algebra of the group
$SL_r$. We choose a Cartan subalgebra $\mathfrak h$ in $\mathfrak g$, and a
Borel subalgebra $\mathfrak b$. This induces a decomposition of $\mathfrak
g$ as $\mathfrak h \oplus \bigoplus_{\alpha \in R} \mathfrak g_{\alpha}$,
where $R$ is the set of \emph{roots} of $\mathfrak g$. The roots appearing
in the decomposition of the Borel subalgebra are said to be
\emph{positive}. 

A choice of an invariant quadratic form on $\mathfrak g$ (a Cartan-Killing
form), allows to define coroots and the lattice of weights $\Lambda$ which
is dual to the lattice generated by coroots. In the case of
$\mathfrak{sl}_r$, roots and coroots are equal: they generate a lattice
which can be naturally identified to 
$\set{(x_1, \dots, x_r) \text{ such that } \sum x_i = 0}$
while the weight lattice is the quotient of $\Z^r$ by the line generated
by the vector $(1, \dots, 1)$.

The positive roots can be chosen to have the form $e_i - e_j$ where $i < j$
and the lattice of weights has a basis of \emph{fundamental weights} 
$\varpi_k = e_1 + \cdots + e_k$ ($e_i$ being the vector having a single
nonzero $i$-th coordinate). 

Let $\lambda = (a_1, \ldots, a_r) \in YD_l$ be a Young diagram. It will
also be considered as a weight of $\mathfrak{sl}_r$. The associated
parabolic subalgebra $\mathfrak p_\lambda$ of $\mathfrak{sl}_r$ is the Lie
algebra spanned by the Borel subalgebra $\mathfrak b$ and the negative root
subspaces corresponding to roots annihilated by $\lambda$. The associated
subgroup of $G := SL_r$ is denoted by $P_\lambda$.

The parabolic weight lattice $\Lambda_\lambda$ is defined to be the
sublattice of $\Lambda$ spanned by the fundamental weights appearing in the
decomposition of $\lambda$: it is the annihilator of the root lattice
annihilated by $\lambda$.

\begin{defn}
  A \emph{parabolic structure} of type $\lambda$ on a $r$-dimensional
  vector space $E$ is the choice of a point in $SL(E)/P_\lambda$.
\end{defn}

\begin{prop}
  If $\lambda^T = (b_1, \ldots, b_l)$, a type $\lambda$ parabolic structure
  on $E$ is equivalent to the datum of a (decreasing) filtration by
  $b_i$-dimensional vector subspaces.
\end{prop}

\begin{proof}
  We have $\lambda = \sum_{i=1}^l \varpi_{b_i}$. Then $\Lambda_\lambda$ is
  the lattice of vectors whose coordinates are constants on blocks of
  length $b_i - b_{i+1}$. It can be shown that $P_\lambda$ is isomorphic to
  the subgroup of $SL_r$ of upper block-triangular matrices where the
  diagonal blocks have size $b_{i+1} - b_i$ exactly. Since $SL(E)$ acts
  transitively on filtrations and $P_\lambda$ is the stabilizer of such a
  structure, one can find a bijection between $G/P_\lambda$ and choices of
  a filtration. 
\end{proof}

The variety $G/P_\lambda$ is projective and carries a natural line bundle
$\mathcal L^\lambda = \mathcal L_{-\lambda}$, obtained as the mixed product
$G \times_{P_\lambda} \C_\lambda$ where $\C_\lambda$ is the one-dimensional
representation of $P_\lambda$ defined by the character $\exp(-\lambda)$.

\begin{prop}
  If $\pi_i$'s are the canonical projections of $G/P_\lambda$ on 
  the Grassmannians $G/P_{\varpi_{b_i}} = \Gr(b_i, E)$, 
  \[ \mathcal L_{-\lambda} = \bigotimes_{i=1}^l \pi_i^\star \O(1) \]
  where the polarisation on the Grassmannians is given by the Plücker
  embedding. Actually
  \[ \vect\pi^\star : \prod_{\text{distinct }i} \Pic(\Gr(b_i, E))
  \to \Pic(G/P_\lambda) \simeq \Lambda_\lambda \]
  is an isomorphism.
\end{prop}

\begin{thm}[Borel-Weil-Bott]
  The vector space of global sections of $\mathcal L_{-\lambda}$ is
  isomorphic as a representation of $G$ to the irreductible representation
  of $G$ of highest weight $\lambda$. 
\end{thm}

By considering $\lambda^T$ instead of $\lambda$, we can define using the
same diagram parabolic structures on a $l$-dimensional vector space $F$. 

\begin{prop}
  \label{parabolic-tensor-product}
  Let $E$ and $F$ be vector spaces of dimensions $r$ and $l$, equipped with
  parabolic structures of type $\lambda$ and $\lambda^T$. Then $E \otimes
  F$ has a natural parabolic structure of type $\varpi_{\abs\lambda}$
  (which we will abusively call a type-$\abs\lambda$ structure).

  This parabolic structure is defined as follows: let
  \[ \begin{aligned}
    0 &\subseteq E_l \subset E_{l-1} \subset\cdots \subset E_1 \subseteq E \\
    0 &\subseteq F_r \subset F_{r-1} \subset\cdots \subset F_1 \subseteq F \\
  \end{aligned} \]
  be the defining filtrations of these structures, then the structure on
  the tensor product is defined by the subspace 
  \[ G = \sum_{i,j} E_i \otimes F_j \]
  where the indices $i,j$ have to satisfy either $j = \lambda^{(i)}$ or $i
  = \lambda^{T,(j)}$. To make it clearer, let $(e_\alpha)$ and $(f_\beta)$ be
  bases of $E$ and $F$ where the chosen Borel subgroups of $SL(E)$ and
  $SL(F)$ are represented by upper triangular matrices. The conditions on
  the indices in the summation can be translated into the condition
  “$e_\alpha \otimes f_\beta$ in the distinguished subspace if and only if
  $\alpha \leq i$ and $\beta \leq j$, where $(i,j)$ are the coordinates of
  a box on the boundary of the diagram of $\lambda$”: in other words, $G$
  is the linear span of the $e_\alpha \otimes f_\beta$ for $(\alpha,\beta)$
  going through the boxes of $\lambda$.

  This construction yields a morphism 
  \[ \tau: SL_r/P_{r,\lambda} \times SL_l/P_{l,\lambda^T} \to 
  SL_{rl}/P_{rl,\abs\lambda} \simeq \Gr(\abs\lambda, rl) \]
  such that $\tau^\star \mathcal L^{\abs\lambda} = \mathcal L^{\lambda}
  \boxtimes \mathcal L^{\lambda^T}$.
\end{prop}

\begin{proof}
  Let $(\mathcal E_i)$ and $(\mathcal F_j)$ be the universal filtrations of
  $\O^r$ and $\O^l$ on $SL_r/P_{r,\lambda}$ and $SL_l/P_{l,\lambda^T}$:
  then $\O^r \boxtimes \O^l$ contains a canonical subbundle $\mathcal G$ of
  rank $\abs\lambda$ defined by the same formula as above, which determines
  a classifying morphism $\tau$ to the Grassmannian of rank $\abs\lambda$
  subspaces of $\C^{rl}$.

  Let $\mathcal G$ denote (abusively) the universal bundle on the
  Grassmannian $\Gr(\abs\lambda, rl)$, which pulls back to $\mathcal G$ on
  $SL_r/P_{r,\lambda} \times SL_l/P_{l,\lambda^T}$. Then $\mathcal
  L^{\abs\lambda}$ is identified with $\det \mathcal G^\vee$ as well as its
  pull-back:
  \[ \tau^\star \mathcal L^{\abs\lambda} = \det \mathcal G^\vee. \]

  The bundle $\mathcal G$ can be filtered by 
  \[ \mathcal G_k = \sum_{i \geq k} E_i \otimes F_{\lambda^{(i)}} \]
  and since
  \[ \mathcal E_{i+1} \otimes \mathcal F_{\lambda^{(i)}} \subset
  \mathcal E_{i+1} \otimes \mathcal F_{\lambda^{(i+1)}} \]
  the $i$-th graded part of $\gr \mathcal G$ is 
  \[ \frac{E_i \otimes F_{\lambda^{(i)}}}
  {(E_i \otimes F_{\lambda^{(i)}}) \cap 
    E_{i+1} \otimes F_{\lambda^{(i+1)}}} = 
  \mathcal E_i / \mathcal E_{i+1} \otimes \mathcal F_{\lambda^{(i)}}. \]
  Now 
  \[ \det \mathcal G = \det \gr \mathcal G = \bigotimes 
  \det(\mathcal E_i / \mathcal E_{i+1})^{\lambda^{(i)}}
  \otimes \det \mathcal F_{\lambda^{(i)}}^{\rank \mathcal E_i - 
    \rank \mathcal E_{i+1}} \]
  and we are reduced to proving this is equal to
  \[ \bigotimes_{i=1}^l \det \mathcal E_i \otimes \bigotimes_{j=1}^r
  \det \mathcal F_j \]
  which can be done counting the number of times $\det \mathcal E_i$ 
  and $\det \mathcal F_j$ appear in each factor. 
\end{proof}

We can give this construction a more functorial look by expressing it in
terms of parabolic maps.

\begin{defn}[Dual parabolic structure]
  Let $E$ be a vector space equipped with a parabolic structure of type
  $\lambda = \sum_{i=1}^l \varpi_{b_i}$ given by a filtration
  $(E_i)_{i=1}^l$. Then its dual $E^\vee$ carries a canonical parabolic
  structure given by the filtration $(E_{l+1-i}^\perp)$ called the \emph{dual
    parabolic structure}, which has type $\lambda^\ast$.
\end{defn}

\begin{defn}[Parabolic map, parabolic morphism]
  Let $E$ be a type $\lambda$ parabolic vector space. Then $\lambda$ as a
  weight of $\mathfrak{sl}_r$ can be written $\sum x_i \varpi_{\alpha_i}$
  for strictly decreasing $0 < \alpha_i < l$ and integer coefficients
  $x_i$. The \emph{reduced filtration} on $E$ is the strictly decreasing
  filtration $(E_i)_{i=1}^s$ defining the parabolic structure.

  Let $F$ be a type $\lambda^{T\ast}$ parabolic vector space. If $\beta_i$
  is $\lambda^{(\alpha_i)}$ then $\lambda^T = \sum y_i \varpi_{b_i}$ with
  distinct increasing $b_i$'s. The reduced filtration $(F_i)_{i=1}^s$ by
  $(l-b_i)$-dimensional subspaces in then \emph{strictly decreasing}. 
  The couples $(\alpha_i,\beta_i)$ correspond to corners of the diagram
  $\lambda$. 

  A \emph{parabolic map} between $E$ and $F$ is a filtered morphism from
  $E$ to $F$, i.e. the image of $E_i$ lies inside $F_i$ for all $i$. We
  denote by $\Hom_{\rm par}(E, F) \subset E^\vee \otimes F$ the vector
  space of parabolic morphisms. 
\end{defn}

\begin{prop}
  \label{parabolic-morphisms}
  Let $E$ be a type $\lambda$ parabolic vector space and $F$ a type
  $\lambda^{T\ast}$ parabolic vector space. The vector space $\Hom_{\rm
    par}(E, F) \subset E^\vee \otimes F$ is precisely the subspace defining
  the type $\abs{\lambda^\ast}$ structure on $E^\vee \otimes F$.
\end{prop}

\begin{proof}
  We consider the reduced filtrations on $E$ and $F^\vee$: since $E$ is
  filtered by $(E_i)$ and $F^\vee$ by $(F_{s+1-i}^\perp)$, the parabolic
  structure on $\Hom(F, E)$ is given by $\sum E_i \otimes F_i^\perp$ (it is
  a type $\abs{\lambda}$ structure), whose annihilator in $\Hom(E, F)$ is
  \[ \bigcap (E_i \otimes F_i^\perp)^\perp
  = \bigcap (E_i^\perp \otimes F + E^\vee \otimes F_i)
  = \bigcap (\Hom(E/E_i, F) + \Hom(E, F_i)) \]
  and each of the spaces being intersected is the subspace of morphisms
  mapping $E_i$ to $F_i$. 
\end{proof}

\subsection{Parabolic vector bundles}

\begin{defn}
  Let $X$ be a smooth curve, and $(\vect p, \vect\lambda)$ a sequence of
  points $p \in D$ with associated Young diagrams $\lambda_p \in YD_l$. A
  \emph{parabolic vector bundle} of type $(\vect p, \vect \lambda)$ is the
  data of a rank $r$ vector bundle $E$ over $X$ and a filtration of
  $E\restriction_p$ by vector subbundles of ranks $\lambda_p^{T,(j)}$
  (i.e. a type $\lambda_p$ parabolic structure on this fiber).

  A parabolic $SL_r$-bundle on $X$ is the data of a principal $SL_r$-bundle
  $E$ on $X$ with the choice of a structure group reduction to
  $P_{\lambda_p}$ (given by a point of $E_p \times_{SL_r} P_{\lambda_p}$)
  at $p$. A parabolic $SL_r$ bundle defined naturally a parabolic
  vector bundle with trivial determinant.
\end{defn}

We know define a notion a section for parabolic bundles which will be
convenient to define the Theta divisors on the corresponding moduli stack. 

\begin{defn}
  If $\vect\lambda$ is a system of Young diagrams, $\O^l$ can be equipped
  with a canonical parabolic structure of type $\vect\lambda^{T\ast}$
  (choosing a filtration by the first coordinates).

  A \emph{parabolic section} of a type $\vect\lambda$ parabolic bundle $E$
  is a parabolic morphism $\O^l \to E(D)$ (i.e. a morphism restricting to a
  parabolic map of vector spaces between the fibers at points $p \in D$). 
  Equivalently, a \emph{parabolic section} is a $l$-tuple of
  sections of $E(D)$ such that the $i$-th section takes values in the $i$-th
  term of the filtration of the fiber at each point. 
\end{defn}

\begin{rem}
  The case $l=1$ and $\lambda=0$ gives back the definition of a
  section of a plain vector bundle. In the general case, parabolic sections
  are allowed to have simple poles along the divisor $D$. 
\end{rem}

\begin{prop}
  The sheaf $E_{\rm par}$ of parabolic sections of $E$ is a
  locally free subsheaf of $E(D)^{\oplus l}$ fitting into a short exact sequence 
  \[ 0 \to E_{\rm par} \to E(D)^{\oplus l} \to 
  \bigoplus_{p \in D,i \in \irange{1,l}} E_p/F_{p,i} \to 0 \]
  where $E_p = E \otimes \O_p$ and $F_{p,j}$ is the $j$-th subspace in the
  filtration of $E_p$ given by its parabolic structure of type $\lambda_p$,
  which is therefore a $\lambda_p^{T,(j)}$-dimensional vector space. The
  dimensions of $E_p/F_{p,i}$ are then the coordinates of
  $\lambda_p^{T,\ast}$.

  The degree of $E_{\rm par}$ is
  \[ l \deg_{\rm par}(E) := l \deg E + \abs{\vect\lambda} \]
  and its slope is $\mu(E_{\rm par}) = \mu(E) + \delta$ where
  $\abs{\vect\lambda} = rl\delta$ and $n$ is the number of marked points.
  The rational number $\mu(E_{\rm par})$ is called the \emph{parabolic
    slope} of $E$, and $\deg_{\rm par}(E) = \deg(E_{\rm par})/l$ its
  \emph{parabolic degree}. 
\end{prop}

\begin{rem}
  Notice that $\abs{\lambda_p} \leq rl$ for any Young diagram of $YD_l$. 
  In particular, $\delta \leq n$. This notion of parabolic degree coincides
  with the usual definition of parablic degree for the parabolic structure
  defined by the weights $(a_1/l, \dots, a_r/l)$ (where $\lambda=(a_1,
  \dots, a_r)$ is the vector corresponding to the Young diagram).
\end{rem} 

\begin{proof}
  The exact sequence is actually a rewriting of the definition and the
  computation of its degree and slope follows from the properties of
  degree with respect to extensions of sheaves.
\end{proof}

\begin{prop}
  If $M$ is a type $\vect\lambda^T$ parabolic vector bundle of rank
  $l$ on $X$, we define the bundle of \emph{parabolic $M$-sections of $E$} to
  be $(M \otimes E)_{\rm par} = \Homrond_{\rm par}(M^\vee, E(D))$. There is an
  exact sequence
  \[ 0 \to (M \otimes E)_{\rm par} \to M \otimes E(D) \to 
  \bigoplus_{p \in D,i \in \irange{1,l}} E_p/F_{p,i} \to 0 \]
  and $\mu((M \otimes E)_{\rm par}) = \mu(E_{\rm par}) + \mu(M)
  = \mu(E) + \mu(M_{\rm par})$.
\end{prop}

\subsection{Theta divisors on moduli spaces of bundles}

Let $\mathcal U_X(r)$ (resp. $\mathcal{SU}_X(r)$) be the moduli stack of
rank $r$ vector bundles with degree zero (resp. trivial determinant) on
$X$. It is an algebraic stack \cite{LMB} and admits a coarse moduli scheme
(which can be seen as the GIT quotient of a Quot scheme) which is a normal
projective variety, whose closed points parameterise $S$-equivalence classes
of semistable vector bundles. It is known \citelist{\cite{DrezetNarasimhan}
  \cite{LaszloSorger}} that the Picard group of the moduli stack (or
equivalently of the coarse moduli scheme, see \cite{BeauvilleLaszlo}) is
isomorphic to $\Z$, as it is the case of the affine Grassmannian
$SL_r(\C((z)))/SL_r(\C[[z]])$.

An ample generator of this group is the \emph{determinant line bundle}
$(\det Rf_\star (\mathcal E \otimes L))^{-1}$ of the “universal” bundle
$\mathcal E$, where $L$ is a degree $g-1$ line bundle on $X$. The resulting
line bundle on $\mathcal{SU}_X(r)$ does not depend on the choice of $L$.

Let $\Upar r$ (resp. $\SUpar r$) be the moduli stack of rank $r$ parabolic
vector bundles with degree zero (resp. trivial determinant), with parabolic
structure of type $\vect\lambda$ at the marked points. There a canonical
forgetful morphism $\phi: \SUpar r \to \mathcal{SU}_X(r)$. Any marked point
$p$ gives rise to a vector bundle $\mathcal E_p$ on $\SUpar r$ which is the
pull-back by the associated section $\SUpar r \to \SUpar r \times \set p
\subset \SUpar r \times X$ of the universal bundle: it is naturally given
a global parabolic structure of type $\lambda_p$.

\begin{prop}[\cite{LaszloSorger}]
  The Picard group of $\SUpar{r}$ is isomorphic to 
  \[ \Pic(\mathcal{SU}_X(r)) \times \prod_{p \in D} \Pic(G/P_{\lambda_p})
  \simeq \Z \times \prod_p \Lambda_{\lambda_p} \]
  where the first factor is generated by the pull-back by $\phi$ of the
  determinant bundle by the fibration $\SUpar{r} \to \mathcal{SU}_X(r)$,
  and the others by the determinants of the subbundles defining the
  filtration of $\mathcal E_p$: these are relative versions of the line
  bundles on $G/P_\lambda$ associated to weights of $P_\lambda$ (which
  are elements of the lattice $\Lambda_\lambda$).
\end{prop}

The moduli stack $\SUpar{r}$ carries a natural ample line bundle $\Lpar l$
represented in the Picard group by the element $(l, \vect \lambda)$ of 
$\Z \times \prod \Lambda_\lambda$.

\begin{prop}
  \label{prop:det-cohomology}
  If $\delta$ is an integer, the natural ample line bundle $\Lpar{l}$ over
  $\SUpar{r}$ is isomorphic to the determinant of the cohomology of $E_{\rm
    par} \otimes L = (E \otimes L)_{\rm par}$ where $L$ is some degree
  $(g-1-\delta)$ line bundle over $X$.
\end{prop}

\begin{proof}
  The proof is similar to \cite{Pauly}: there is an exact sequence 
  \[ 0 \to (\mathcal E \otimes L)_{\rm par} \to (\mathcal E \otimes L(D))^l
  \to \bigoplus_{p \in D,i \in \irange{1,l}} E_p/F_{p,i} \to 0 \]
  and if $\pi$ is the projection from $\SUpar{r} \times X$ to $\SUpar{r}$
  we have
  \[ \det R\pi_\star(\mathcal E_{\rm par} \otimes L)^\vee
  = \det R\pi_\star(\mathcal E \otimes L(D))^{\otimes -l}
  \otimes \bigotimes_{p \in D} \bigotimes_{i=1}^l
  \det(\mathcal E_p/\mathcal F_{p,i}) \]
  and $\det(\mathcal E_p/\mathcal F_{p,i}) \simeq \det \mathcal
  F_{p,i}^\vee$ are exactly the positive generators of the
  $\Lambda_{\lambda_p}$ component of the Picard group as above.
\end{proof}

In the rest of this section, we suppose $\delta$ is an integer.

\begin{prop}
  For any type $\vect\lambda^T$ parabolic vector bundle $M$ of rank $l$ and
  parabolic slope $g-1$, the determinant of the cohomology of $(M \otimes
  \mathcal E)_{\rm par}$ is a line bundle isomorphic to $\Lpar l$ and
  contains a canonical section whose vanishing detects the existence of
  parabolic $M$-sections.
\end{prop}

\begin{rem}
  The previous proposition is simply the particular case where $M =
  L^{\oplus l}$ is equipped with the type $\vect\lambda^T$ parabolic
  structure of $\C^{l}$.
\end{rem}

\begin{cor}
  \label{prop:canonical-section}
  Let $M$ be a fixed parabolic bundle of type $\vect\lambda^T$. There is a
  canonical section (up to homothety) in $H^0(\Lpar l)$ whose vanishing
  locus on the moduli stack consists of of parabolic bundles having nonzero
  parabolic sections (after tensoring by $L$).
\end{cor}

\begin{proof}
  If $\pi$ is the projection from $\SUpar{r} \times X$ to $\SUpar{r}$,
  $\Lpar{l} = \det R\pi_\star (\mathcal E \otimes L)_{\rm par}^\vee$ 
  can be represented by the determinant of any quasi-isomorphic complex of
  vector bundles with two terms: it has a canonical section given by the
  determinant of the differential, whose vanishing detects the existence of
  sections (which are the cohomology of this complex). 
\end{proof}

\begin{prop}
  When $l = 1$, the canonical section associated to a line bundle $L$ of
  degree $g-1-\delta$ is nonzero and defines a divisor on the moduli
  stack.
\end{prop}

\begin{proof}
  When $l=1$, a parabolic structure of level $1$ on a vector bundle is the
  choice of a single subspace at each of the marked point. A parabolic
  bundle $E \otimes L$ has nonzero parabolic sections if and only if
  $E_{\rm par} \otimes L$ has nonzero ordinary sections (here $E_{\rm par}$
  is a subvector bundle of $E(D)$. Let $F$ be a vector bundle of degree
  $\delta$ such that $F \otimes L$ has no nonzero section. Such a bundle
  exists because the Theta map at level 1 is defined on an open set. Then
  for any extension $E$ of $F$, with a parabolic structure such that
  $E_{\rm par} = F$, the canonical section defined by $L$ does not vanish
  at $E$.
\end{proof}

We will often refer to the “canonical section” without mentioning that it
is only defined up to a scalar.

The determinant line bundle can also be defined on $\Upar{r}$ but it
depends on the choice of $M$. Let $M$ be a type $\vect\lambda^{T\ast}$
parabolic bundle of parabolic slope $g-1$ over $X$ and rank $l$. We set
$\Lpar{l}_M = \det R\pi_\star (M \otimes \mathcal E)_{\rm par}^\vee$.

\begin{prop}
  The line bundle $\Lpar{l}_M$ is isomorphic to
  \[ \phi^\star \Theta^l_{r,M} \otimes \Lpar{0} \otimes \bigotimes_p 
  (\det \mathcal E_p)^{\otimes l} \]
  where $\Theta_{r,M}$ is the determinant line bundle defined by $M$ on
  $\mathcal U_X(r)$. 
\end{prop}

\begin{proof}
  The proof is identical to \ref{prop:det-cohomology}, except $\mathcal
  E_p$ no longer has a trivial determinant. 
\end{proof}

\begin{prop}
  The line bundles $\Lpar{l}_M$ satisfy the following properties:
  if $M$ and $N$ are vector bundles of parabolic slope $g-1$, then 
  \begin{itemize}
  \item $\mathcal L^{sl,s \vect\lambda}_{M^{\oplus s}} = (\Lpar{l}_M)^s$
  \item $\Lpar{l}_M = \Lpar{l}_N \otimes 
    {\det}^\star(\det M \otimes \det N^{-1})$.
  \end{itemize}
  Here a line bundle of degree zero acts by translation on the Picard
  varieties $\Jac^d(X)$ and its action on line bundles is isomorphic to the
  tensor product by some line bundle on the Picard variety, which will be
  denoted similarly. 
\end{prop}

\begin{proof}
  The first equality results from the definitions. The second one is a
  consequence of the identity
  \[ \Theta_{r,M} =  \Theta^l_{r,L} \otimes {\det}^\star(\det M \otimes
  \det L^{-l}) \]
  whree $L$ is a line bundle of the same slope as $M$. 
  It will be often used to compute decompositions of pull-backs
  of determinant line bundles.
\end{proof}

\begin{prop}
  The tensor product of parabolic structures extends to parabolic vector
  bundles, giving a morphism of stacks 
  \[ \tau: \SUpar{r} \times \Utpar{l} \to \Uppar{rl} \]   
  mapping a couple of parabolic vector bundles $(E, F)$ to their 
  parabolic tensor product $E \otimes F$.

  Furthermore, for any line bundle $L$ of degree $g-1-\delta$
  \[ \tau^\star \Lppar{1}_L = \Lpar{l} \boxtimes \Ltpar{r}_L. \]

  The pull-back of the canonical section of $\Lppar{1}$ defines 
  a pairing $SD$ between $H^0(\SUpar{r}, \Lpar{l})^\vee$
  and $H^0(\Utpar{l}, \Lpar{r}_L)^\vee$ by duality.
\end{prop}

\begin{proof}
  The stack $\SUpar{r} \times \Utpar{l}$ carries a parabolic vector bundle
  $\mathcal E \boxtimes \mathcal F$ where $\mathcal E$ and $\mathcal F$ are
  the tautological universal vector bundles: since this is a type
  $\abs{\vect\lambda}$ and rank $rl$ parabolic vector bundle, there is an
  associated classifying morphism to $\Uppar{rl}$.

  The pull-back formula can be derived for the formula in the non-parabolic
  case (which has the form $\tau^\star \mathcal L_{rl,L} = \mathcal L_r^k
  \boxtimes \mathcal L_{k,L}^r$) and from the computation of the pull-back of 
  $\mathcal L^{0,\vect{\abs\lambda}}$ which is identical to the case of
  vector spaces (proposition \ref{parabolic-tensor-product}).
\end{proof}

The pairing $SD$ is the classical strange duality map. Following Marian and
Oprea \cite{MarianOprea}, we now introduce new spaces of sections which
satisfy more symmetries and will benefit from the enumerative geometry
constructions. 

\begin{prop}
  \label{prop:sd2-definition}
  There is a morphism
  \[ (\tau,\delta): \Upar{r} \times \Utpar{l} \to \Uppar{rl} \times \Jac
  X \]
  mapping a couple of parabolic vector bundles $(E, F)$ to 
  $(E \otimes F, \det E^\vee \otimes \det F)$.

  Furthermore, if $\Theta$ is a Theta divisor on $\Jac X$, and $\bar\Theta$
  its image by the involution $x \mapsto -x$, a choice of a degree
  $g-1-\delta$ line bundle $L$ gives lines bundles on the moduli stacks
  satisfying:
  \[ (\tau,\delta)^\star (\Lppar{1}_L \boxtimes \Theta)
  = (\Lpar{l}_L \otimes {\det}^\star \bar\Theta) \boxtimes 
  (\Ltpar{r}_L \otimes {\det}^\star \Theta). \]

  The pull-back of the external product of canonical sections of 
  $\Lppar{1}_L$ and $\Theta$ defines a pairing $\widetilde{SD}$ between 
  $H^0(\Upar{r}, \Lpar{l}_L \otimes {\det}^\star \bar\Theta)^\vee$ and 
  $H^0(\Utpar{l}, \Ltpar{r}_L \otimes {\det}^\star \Theta)^\vee$.
\end{prop}

\begin{proof}
  The decomposition of the pull-back of the determinant bundle is an
  application of Mumford's see-saw principle, and is analogous to the
  non-parabolic case \cite{MarianOprea}: let $E$ be a closed point of
  $\Upar{r}$. The restriction of $(\tau,\delta)^\star (\Lppar{1}_L
  \boxtimes \Theta)$ to the fiber $\set{E} \times  \Utpar{l}$ is 
  \[ \det R\pi_\star(E \otimes \mathcal F \otimes L)_{\rm par}
  \otimes {\det}^\star (t_{\det E^\vee}^\star \Theta) \]
  where $t_\bullet$ is the translation operator. This is equal to
  \[ \det R\pi_\star(\mathcal F \otimes L)_{\rm par}
  \otimes {\det}^\star(\det E) \otimes {\det}^\star (\det E^\vee \otimes
  \Theta) = \Ltpar{r}_L \otimes {\det}^\star \Theta. \]
  
  A similar formula holds for $\Upar{r} \times \set{F}$ fibers.
\end{proof}

\begin{prop}
  If $\widetilde{SD}$ is a perfect pairing for some choice of $\Theta$,
  then so is $SD$.
\end{prop}

\begin{proof}
  Let $\theta$ be the unique (up to homothety) section of $\bar\Theta$ on
  $\Jac X$, also denoting its pullback to $\Utpar{l}$, and denoting quite
  abusively the multiplication morphism $\Ltpar{r}_L \to \Ltpar{r}_L \otimes
  {\det}^\star \bar\Theta$.
  
  Since $\widetilde{SD}$ is induced by a section $\Delta \cdot \theta$
  while $SD$ is only induced by $\Delta$, we have the following commutative
  diagram:
  \[ \xymatrix{
    H^0(\SUpar{r}, \Lpar{l})^\vee \ar[r]^{\rho^\dagger}
    \ar[d]^{SD} & H^0(\Utpar{l}, \Ltpar{r}_L) \ar[d]^{1 \otimes \theta} \\
    H^0(\Upar{r}, \Lpar{l} \otimes {\det}^\star \Theta)^\vee
    \ar[r]^{\widetilde{SD}}
    & H^0(\Utpar{l}, \Ltpar{r}_L \otimes {\det}^\star \bar\Theta) \\
  } \]
  \noindent where $\rho$ is the restriction from $\Upar r$ to $\SUpar r$, and
  $\rho^\dagger$ its transpose. Since $SD$ et $\widetilde{SD}$ are maps
  between spaces of the same dimension, as we will see in section
  \ref{section:rank-lvl-symmetry}, if $\rho$ is surjective, $\rho^\dagger$
  is injective and the conclusion follows.

  The existence of the étale covering $\SUpar{r} \times \Jac X \to
  \Upar{r}$, allows to extend a section of $\Lpar l$ from $\SUpar{r}$ to
  $\Upar l$, by tensoring by a section of a theta bundle and projecting to
  the space of invariants under translation by the $r$-torsion subgroup, so
  that it descends to $\Upar{r}$. This averaging technique works for a
  generic choice of the chosen theta section \cite[Proposition
  4]{MarianOprea}.
\end{proof}

\subsection{Moduli spaces of non-zero degree bundles}
\label{subsection:arbitrary-degree}

Let $d$ be an integer and $\Dcyr$ be a line bundle of degree $d$. Then we
can define the moduli stack of degree $d$ vector bundles on $X$, $\mathcal
U_X(r, d)$ and its parabolic counterpart $\Upar{r,d}$, for a choice of
Young diagrams $\vect\lambda$. We define similarly moduli stacks
$\SUpar{r,\Dcyr}$ of vector bundles with fixed determinant $\Dcyr$. These
are (non-canonically) isomorphic to each other, and $\Upar{r,d}$ is
isomorphic to $\Upar{r,d+r}$ by the map $E \to E \otimes \Lcyr$ for 
any degree 1 line bundle $\Lcyr$.

The formula for the parabolic slope shows that we can still define a
determinant line bundle $\Lpar{l}_M$ on $\Upar{r,d}$ with a canonical
section if $\mu_{\rm par} = d/r + \abs{\vect\lambda}/(rl) \in (1/l)\Z$
in which case we choose a parabolic vector bundle $M$ of rank $l$ and
ordinary slope $(1-g+\mu_{\rm par})$, so that $\Homrond_{\rm par}(M^\vee,E(D))$
has slope $g-1$ for any $E \in \Upar{r,d}$. This condition is equivalent to 
\[ \abs{\vect\lambda} + ld \equiv 0 \mod r \text{ or }
\abs{\vect\lambda^\ast} \equiv ld \mod r. \]

\begin{prop}
  Let $d$ and $\dcyr$ be integers, and $\vect\lambda$ be a system of Young
  diagrams such that there exists an integer $\mu$:
  \[ \abs{\vect\lambda} + ld+r\dcyr = rl\mu. \]
  Then $\abs{\vect\lambda^\ast}$ is congruent to $ld\textup{ mod }r$ and to
  $r\dcyr\textup{ mod }l$. We define a tensor product morphism
  \[ \widetilde\tau: \Upar{r,d} \times \Utpar{l,\dcyr} \to
  \Uppar{rl, ld+r\dcyr} \times \Pic^{\dcyr-d} X \]
  mapping $(E,F)$ to their parabolic tensor product $(E \otimes F, \det
  E^\vee \otimes \det F)$. Let $\Lcyr$ be a line bundle with slope
  $g-1-\mu$ and $\Theta$ be a Theta divisor on $\Pic^{\dcyr-d} X$. Then
  choosing a base point $(E_0, F_0)$ defines a divisor $\overline\Theta_d$
  on $\Pic^d X$ pulling back by $x \to \det F_0 - x$ and $\Theta_\dcyr$ on
  $\Pic^\dcyr X$ pulling back by the action of $\det E_0^\vee$.

  Then
  \[ \widetilde\tau^\star \Lppar{1}_\Lcyr \boxtimes \Theta =
  (\Lpar{l}_{F_0 \otimes \Lcyr} \otimes {\det}^\star \overline\Theta_d)
  \boxtimes
  (\Ltpar{r}_{E_0 \otimes \Lcyr} \otimes {\det}^\star \Theta_\dcyr) \]

  Moreover there exists a canonical section in $H^0(\Uppar{rl, ld+r\dcyr},
  \Lppar{1}_\Lcyr) \otimes {\det}^\star \Theta$ whose pull-back defines a
  pairing $\widetilde{SD}$ between $H^0(\Upar{r,d}, \Lpar{r}_{F_0 \otimes
    \Lcyr} \otimes {\det}^\star \overline\Theta_d)^\vee$ and 
  $H^0(\Utpar{l,\dcyr}, \Ltpar{l}_{E_0 \otimes \Lcyr} \otimes 
  {\det}^\star \Theta_\dcyr)^\vee$.
\end{prop}

\begin{rem}
  In the case $d = \dcyr = 0$, choosing $E_0 = \O_X^r$ and $F_0=\O_X^l$
  gives proposition \ref{prop:sd2-definition}.
\end{rem}

\begin{proof}
  The proof follows the same steps as proposition
  \ref{prop:sd2-definition}. Let $E$ be a closed point of $\Upar{r,d}$. The
  restriction of $\widetilde\tau^\star (\Lppar{1}_\Lcyr \boxtimes \Theta)$
  to the fiber $\set{E} \times \Utpar{l,\dcyr}$ is
  \[ \det R\pi_\star(E \otimes \mathcal F \otimes \Lcyr)_{\rm par}
  \otimes {\det}^\star (t_{\det E^\vee}^\star \Theta) \]
  where $t_\bullet$ is the translation operator. We decompose 
  $\det E^\vee = \det E_0^\vee \otimes (\det E_0 \otimes \det E^\vee)$ and
  get 
  \[ \det R\pi_\star(E \otimes \mathcal F \otimes \Lcyr)_{\rm par} \otimes
  {\det}^\star (\Theta_\dcyr \otimes  (\det E_0 \otimes \det E^\vee)) \]
  which is equal to
  \[ \det R\pi_\star(E_0 \otimes \mathcal F \otimes \Lcyr)_{\rm par}
  \otimes {\det}^\star \Theta_\dcyr. \]
  
  The proof for fibers of type $\Upar{r,d} \times \set{F}$ is identical.
\end{proof}

The proof of the following proposition follows the degree zero case. 

\begin{prop}
  Let $\vect\lambda$ be a system of Young diagrams satisfying the same
  conditions. The tensor product of parabolic bundles defines a morphism
  \[ \tau: \Upar{r,d} \times \SUtpar{l,\Dcyr} \to 
  \Uppar{rl,ld+r\dcyr} \]   
  where $\dcyr = \deg \Dcyr$. If $F_0$ is a given point in
  $\SUtpar{l,\Dcyr}$ and $\Lcyr$ a line bundle with appropriate slope, we
  have
  \[ \tau^\star \Lppar{1}_\Lcyr = \Lpar{l}_{F_0 \otimes \Lcyr} \boxtimes 
  \Ltpar{r} \]
  which does not depend on $F_0$ (since it always has determinant
  $\Dcyr$).

  The pull-back of the canonical section of $\Lppar{1}$ defines 
  a pairing $SD$ between $H^0(\Upar{r}, \Lpar{l}_{F_0 \otimes \Lcyr})^\vee$
  and $H^0(\SUtpar{l}, \Lpar{r})^\vee$ by duality.
\end{prop}

As in the degree zero case, if $\widetilde{SD}$ is an isomorphism for a
suitable choice of $\Theta$, then $SD$ is an isomorphism too. We still have
to check the dimensions are equal, which will be done in section
\ref{section:symmetries}.

\subsection{Parabolic structures and Schubert varieties}

In this section, we fix a $N$-dimensional vector space $\Gamma$ where
$N=r+l$, and denote by $B \subset SL(\Gamma)$ a choice of Borel subgroup in
$SL(\Gamma)$ (which is the stabiliser of a given complete flag $(\Gamma_i)$
in $\Gamma$).

Let $\lambda$ be a Young diagram in a $l \times r$ rectangle. The following
notations will be used:
\[ \begin{aligned}
  \lambda &= (a_1^{b_k} a_2^{b_{k-1}-b_k} \cdots a_k^{b_1-b_2}) \\
  I_\lambda &= (l+1, \cdots, l+r) - \lambda \\
  \lambda^T = \mu &= (b_1^{a_k} b_2^{a_{k-1}-a_k} \cdots
  b_k^{a_1-a_2}) \\
  J_\mu &= (l, \cdots, 1) + \mu \\
\end{aligned} \]
where $(a_i)$ and $(b_j)$ are decreasing sequences of integers, which are
the lengths of rows and columns of $\lambda$. With these notations,
$\lambda \mapsto I_\lambda$ is a bijection between Young diagrams less than
$(l^r)$ and $r$-element subsets of $\irange{1,r+l}$, and $I_\lambda$ and
$J_\mu$ are complements of each other.

A subset $I$ of $\irange{1,r+l}$ define a partial flag $(\Gamma_i)_{i \in
  I}$ of $\Gamma$.

\begin{defn}[\cite{Fulton-schubert}]
  The \emph{Schubert variety} $\overline{Y_\lambda}$ is the codimension
  $\abs\lambda$ subvariety of $\Gr(r, \Gamma)$ whose points are the
  $r$-dimensional subspaces $E$ of $\Gamma$ such that the inclusion of $E$
  in $\Gamma$ is a morphism of filtered vector spaces for the filtration
  given by some complete flag of $E$ and the partial flag $(\Gamma_i)_{i
    \in I_\lambda}$.

  The \emph{Schubert cell} $Y_\lambda$ is the open set of
  $\overline{Y_\lambda}$ where there is a unique complete flag of $E$
  compatible with the inclusion as above. It is isomorphic to an affine
  space.
\end{defn}

Let $(e_i)$ be a triangular basis of $\Gamma$ with respect to $B$ (the
subspaces spanned by the first vectors of the basis are stable under $B$),
and $(e'_i)$ its dual basis (with reversed numbering so that $\pair{e'_i,
  e_j} = \delta_{i+j,N+1}$), which is triangular with respect to the
negative Borel subgroup $\bar B$ in $GL(\Gamma^\vee) \simeq GL(\Gamma)$.
We consider and $\Gamma^\vee$ the corresponding complete flag:
\[ 0 = \Gamma'_0 = \Gamma_N^\perp \subset \cdots \subset
\Gamma'_N = \Gamma_0^\perp = \Gamma^\vee. \]

If $T$ denotes the diagonal maximal torus corresponding to the basis
$(e_i)$, the $T$-fixed points in $\Gr(r, \Gamma)$ are the $r$-dimensional
subspaces spanned by a subset of the basis: they can be indexed by the
$I_\lambda$ so we will write $e_{I_\lambda}$ or $e_\lambda$ for these
subspaces. Similarly, $e'_\mu$ denote the $T$-fixed points of
$\Gr(k,\gamma^\vee)$.

We have $Y_\lambda = Be_\lambda$ and $Y'_\mu = \bar B e'_\mu$ for Schubert
varieties of $\Gamma$ and $\Gamma^\vee$. Mapping a subspace to its
annihilator yields an isomorphism of $Y_\lambda \subset Gr(r, \Gamma)$ onto
$Y'_\mu \subset Gr(l, \Gamma^\vee)$.

\begin{prop}
  The tautological bundle $\mathcal E \subset \Gamma \otimes \O$ over
  $Y_\lambda$ has a natural filtration which makes it a $P_\lambda$-bundle.
\end{prop}

\begin{proof}
  It is enough to find rank $b_j$ subbundles. But by the very definition of
  the Schubert cell, checking the $b_j$-th Schubert condition gives 
  precisely that $\mathcal E \cap \Gamma_{l+b_j-a_{k+1-j}}$ is at every point a
  $b_j$-dimensional subspace of the fiber of $\mathcal E$.
\end{proof}

\begin{rem}
  Considering the torus-fixed point $e_{I_\lambda}$ gives a clear view of
  this parabolic structure: it corresponds to the decomposition of
  $I_\lambda$ in sequences of consecutive integers. 
\end{rem}

\begin{prop}
  The tautological exact sequence 
  \[ 0 \to \mathcal E \hookrightarrow 
  \Gamma \otimes \O \to \mathcal F \to 0 \]
  has a natural structure of an exact sequence of \emph{parabolic vector 
    bundles}, meaning that it is an exact complex of filtered vector
  bundles, the filtration being induced by the $P_\lambda$-structure of
  $\mathcal E$, the $P_{\lambda^{T\ast}}$-structure of $\mathcal F$ and the
    canonical filtration of $\Gamma$ by $(a_i^\ast+b_i)$-dimensional subspaces
    (given by the boundary of the Young diagram), where $a_i^\ast =
    l-a_{k+1-i}$ are the row lengths of $\lambda^\ast$.
\end{prop}

\begin{proof}
  Recall that the filtration on $\mathcal E$ is induced by 
  $\Gamma_{l+b_i-a_{k+1-i}} = \Gamma_{a_i^\ast+b_i}$. Considering the dual 
  sequence 
  \[ 0 \to \mathcal F^\vee = \mathcal E^\perp \hookrightarrow
  \Gamma^\vee \otimes \O \to \mathcal E^\vee \to 0 \]
  gives the answer for the second arrow. 
\end{proof}

\begin{prop}
  Let $0 \to \mathcal E \to \Gamma \otimes \O \to \mathcal F \to 0$ be the
  tautological exact sequence on $Y_\lambda$. Let $\mathcal E^\perp =
  \mathcal F^\vee$ be the annihilator of $\mathcal E$ in $\Gamma^\vee
  \otimes \O$ (it is the tautological subbundle over $Y'_\mu$).

  Then the conormal bundle of $Y_\lambda$ in $\Gr(r, \gamma) \simeq \Gr(l,
  \Gamma^\vee)$ is naturally isomorphic to the canonical subbundle of 
  $\Homrond(\mathcal F, \mathcal E) = \mathcal E \otimes \mathcal E^\perp$
  induced by their parabolic structures of type $\lambda$ and
  $\lambda^T$. 

  Consequently the tangent bundle to $Y_\lambda$ is the subbundle
  $\Homrond_{\rm par}(\mathcal E, \mathcal F)$ of $\Homrond(\mathcal E,
  \mathcal F)$ (which is the tangent bundle to the Grassmannian) consisting
  of parabolic morphisms.
\end{prop}

\begin{proof}
  Equivariance under the action of the Borel subgroup $B \subset
  SL(\gamma)$ can be used to prove the result once it is known at the
  torus-fixed point $E = e_I = e_{I_\lambda}$ and its annihilator 
  $F^\vee = E^\perp = e'_J$.

  Using \ref{parabolic-morphisms} it is equivalent to prove the statement
  concerning the conormal bundle and the one about the tangent bundle.

  We can parameterise $Y_\lambda$ by bases adapted to the unique full flag
  satisfying the Schubert conditions: they can be represented naturally as
  block-triangular rectangular matrices, and this defines an isomorphism of
  $T_E Y_\lambda$ with $Y_\lambda$. It is easy to see that these triangular
  matrices represent exactly the parabolic morphisms.
\end{proof}

\section{Computation of the Verlinde numbers}

\subsection{Verlinde numbers and the Verlinde formula}

The conformal blocks are finite-dimensional vector spaces combining
representation theory of affine Lie algebras and the geometry of $X$. They
were proved to be isomorphic to the spaces of sections of line bundles on
the moduli space of (parabolic) vector bundles
\citelist{\cite{BeauvilleLaszlo} \cite{Faltings} \cite{Pauly-blocs}}.

Their dimension is computed by the Verlinde formula \cite{Verlinde} which
can be used in a much more general case to obtain the dimensions of spaces
involved in Rational Conformal Field Theories. Its mathematical proof
relies upon the degeneration properties proved in \cite{TUY}, and more
explicit formulae can be obtained by representation-theoretic arguments
\cite{Beauville-verlinde}. Computational details about the combinatorial
properties of the Verlinde numbers can be found in \cite{Zagier} for the
$SL_r$ case and \cite{OxburyWilson} for a more general case.

\begin{prop}
  \label{Verlinde-formula}
  Assume $\abs{\vect\lambda}$ is divisible by $rl$. Then the dimension of the
  space of parabolic conformal blocks (sections of the natural ample line
  bundle defined by $\vect\lambda$) is
  \[ \begin{aligned}
    h^0(\SUpar{r}, \Lpar{l}) &= \frac{r^g}{(r+l)^g} \\
    &\times \sum_{S \in \mathfrak P_r(\Z/(r+l)\Z)} 
    \prod_{s \in S, t \notin S} \abs{2 \sin \pi\frac{s-t}{r+l}}^{g-1}
    \prod_{s \in S} (\zeta^s)^{-\abs{\vect\lambda}/r} 
    \prod_p S_{\lambda_p}(\zeta^s, s\in S)
  \end{aligned} \]
  where $\zeta$ is a primitive $(r+l)$-th root of unity and $\mathfrak
  P_r(\Sigma)$ denote the set of $r$-element subsets of $\Sigma$.

  The following formula also holds:
  \[ \begin{aligned}
    h^0(\Upar{r}, \Lpar{l}_M) &= \frac{l^g}{(r+l)^g} \\
    &\times \sum_{S \in \mathfrak P_r(\Z/(r+l)\Z)} 
    \prod_{s \in S, t \notin S} \abs{2 \sin \pi\frac{s-t}{r+l}}^{g-1}
    \prod_{s \in S} (\zeta^s)^{-\abs{\vect\lambda}/r} 
    \prod_p S_{\lambda_p}(\zeta^s, s\in S)
  \end{aligned}\]
  for sections over the moduli stack of degree zero parabolic vector
  bundles, and any choice of a parabolic vector bundle $M$ of type
  $\vect\lambda^T$.
\end{prop}

\begin{proof}
  We use the notations of \cite{Beauville-verlinde} where the Verlinde
  formula in the particular case of $SL_r$ can be written
  \[ h^0(\SUpar{r}, \Lpar{l}) = \abs{T_l}^{g-1} 
  \sum_{t \in T_l^{\mathrm{reg}}/W}
  \frac{\prod_p \Tr_{V_{\lambda_p}}(t)}{\Delta(t)^{g-1}}. \]
  We choose a maximal torus $T$ and a Borel subgroup of $SL_r$, and 
  the Weyl character formula provides $\Tr_{V_\lambda}(t) =
  S_\lambda(t)$ (for $t \in T$ considered as the $r$-tuple of its
  eigenvalues) and $\Delta(t) = \abs{\Vdm(t)}^2$.

  The set $T_l$ is the subgroup of the maximal torus such that $t^{r+l}$ is
  a scalar matrix: it is the subgroup annihilated by the sublattice
  $(l+h^\vee)Q$ where $Q$ is the weight lattice and $h^\vee = r$ is the
  dual Coxeter number of $A_{r-1}$. Its order is
  $\frac{r}{r+l}(r+l)^r$. The set $T_l^{\mathrm{reg}}$ is the subset of
  $T_l$ of elements which are not stabilised by any element of the Weyl
  group. Let $\overline{T_l}^{\mathrm{reg}}$ be its image in $PSL_r$.

  Let $\bar t$ be an orbit in $\overline{T_l}^{\mathrm{reg}}/W$. It
  corresponds to:
  \begin{itemize}
  \item $r$ elements of $T_l^{\mathrm{reg}}$, related by multiplication by 
    $r$-th roots of unity;
  \item $(r+l)$ subsets of size $r$ in $\Z/(r+l)\Z$, which can be 
    mapped to elements of $GL_r$ by $\exp\pa{\frac{2i\pi}{r+l} \bullet}$.
  \end{itemize}

  Since $\prod S_\lambda$ is a homogeneous symmetric polynomial of total
  degree $\abs{\vect\lambda}$ being a multiple of $rl$, its value at $t$ 
  is unchanged when $t$ is multiplied by a $r$-th root of unity: it only
  depends on the equivalence class $\bar t$.

  A diagonal matrix $\zeta^S \in GL_r$ can be replaced by an element $t$ of
  $SL_r$ multiplying by a $r$-th root of its determinant. The value of
  $S_{\vect\lambda}$ does not depend of this choice, so if $S$ is a subset
  of $\Z/(r+l)\Z$, 
  \[ S_{\vect\lambda}(\zeta^S) = S_{\vect\lambda}(t) 
  \pa{\prod \zeta^s}^{-\abs{\vect\lambda}/r}, \]
  giving the formula
  \[ \begin{aligned}
    h^0(\SUpar{r}, \Lpar{l}) &= \pa{\frac{r}{r+l}}^{g-1} 
    (r+l)^{r(g-1)} \frac{r}{r+l} \\
    &\times\sum_{S \in \mathfrak P_r(\Z/(r+l)\Z)} 
    \pa{\prod_{s \in S} \zeta^s}^{-\abs{\vect\lambda}/r}
    \prod_i S_{\lambda_i}(\zeta^S) \abs{\Vdm(t)}^{2-2g}.
  \end{aligned} \]

  Given a fixed subset $S$, the following holds:
  \[ \abs{\Vdm(\zeta^s, s \in S)}^2 
  = \prod_{s \neq s' \in S} \abs{\zeta^s - \zeta^{s'}} 
  = \prod_{s \in S} \frac{r+l}{\prod_{s' \notin S} 
    \abs{\zeta^s - \zeta^{s'}}}
  = \frac{(r+l)^r}{\prod_{s \in S, s' \notin S} \abs{\zeta^s - \zeta^{s'}}} \]
  and writing $\zeta = \exp\pa{\frac{2i\pi}{r+l}}$, $\abs{\zeta^s - \zeta^t}
  = 2 \sin \pi \frac{s-t}{r+l}$, we recover the desired formula. 

  The second equality can be derived from the étale cover of degree
  $r^{2g}$ \cite{DonagiTu}, 
  \[\tau: \SUpar{r} \times \Jac(X) \to \Upar{r}\]
  corresponding to the tensor product functor. Mumford's see-saw principle
  shows that $\tau^\star \Lpar{l}_M = \Lpar{l} \boxtimes \Theta^{rl}$ for
  some theta divisor $\Theta$ (any parabolic bundle in $\SUpar{r}$ defines
  the same determinant bundle on $\Jac X$) hence
  \[ r^{2g} h^0(\Upar{r}, \Lpar{l}_M) = (rl)^g h^0(\SUpar{r}, \Lpar{l}). \]
\end{proof}

An interesting consequence of this formula is that the theta divisor is the
only element of the linear system given by $\Lpar 1$ if $l = 1$. 

\begin{cor}
  If $l = 1$ and $\vect\lambda$ only consists of level $1$ weights
  (fundamental weights), the vector space $H^0(\Upar{r}, \Lppar{1})$ has
  dimension one if $\abs{\vect\lambda}$ is a multiple of $r$.
\end{cor}

\begin{proof}
  The level $1$ Schur polynomials are actually the elementary symmetric
  functions. The sum in the Verlinde formula involves the subsets $S$ whose
  complement contains exactly one element: now if $S = \Z/(r+1)\Z
  \setminus \set{s}$ the values of the elementary symmetric functions on
  $\zeta^S$ are the coefficients of the polynomial
  \[ \frac{X^{r+1} - 1}{X - \zeta^s} = 
  \sum_{k = 0}^{r} (-1)^k S_k(\zeta^S) X^{r-k} \]
  thus $S_k(\zeta^S) = (-1)^k \zeta^{sk}$ (here we denote by $S_k$ the
  polynomial $S_{\varpi_k}$). We hence get
  \[ \prod_{\abs\lambda} S_{\abs\lambda}(\zeta^S) 
  = (-\zeta^s)^{\abs{\vect\lambda}}. \]
  The other factor is
  \[ \prod_{\sigma \in S} (\zeta^\sigma)^{-\abs{\vect\lambda}/r}
  = ((-1)^r \zeta^{-s})^{-\abs{\vect\lambda}/r}
  = (-1)^{\abs{\vect\lambda}} (\zeta^s)^{\abs{\vect\lambda}/r} \]
  because $\prod_{s=0}^{r} \zeta^s = (-1)^r$.
  And since $\abs{\vect\lambda} + \abs{\vect\lambda}/r$ is a
  multiple of $(r+1)$, these two factors cancel each other. 
  
  The only remaining factor is the one which does not depend on
  $\vect\lambda$, and we have actually 
  \[ h^0(\Upar{r}, \Lpar{1}) = h^0(\mathcal U_X(r), \mathcal L) \]
  where $U_X(r)$ is the moduli stack of plain vector bundles of rank $r$
  and degree $0$ on $X$ and this number equals
  \[ \begin{aligned}
    h^0(\mathcal U_X(r), \mathcal L) &= 
    \frac{1}{(r+1)^g} \sum_{s = 0}^r \prod_{t \neq s} 
    \abs{2 \sin \pi\frac{s-t}{r+1}}^{g-1} \\
    &= \frac{1}{(r+1)^{g-1}} \prod_{t = 1}^r
    \abs{2 \sin\pi\frac{t}{r+1}}^{g-1} = 1.
  \end{aligned} \]
\end{proof}

\subsection{Rank-level symmetry between Verlinde numbers}
\label{section:rank-lvl-symmetry}

The rank-level symmetry is a consequence of the well-known reciprocity law
between the expressions of Schur polynomials in terms of complete symmetric
polynomials and elementary symmetric polynomials \cite{Macdonald}. 

\begin{lemma}
  Let $S$ be a subset consisting of $(r+l)$-th roots of unity with $r$
  elements, and $T$ be its complement. Then
  \[ S_\lambda(u, u \in S)
  = (-1)^{\abs\lambda} S_{{\lambda^\ast}^T}(\bar v, v \in T)
  = (-1)^{\abs\lambda} S_{\lambda^T}(v, v \in T). \]
\end{lemma}

\begin{proof}
  Let $\mathcal F$ denote the matrix with coefficients
  $(\zeta^{ij}/\sqrt{r+l})_{i,j}$, where $\zeta$ is a primitive $(r+l)$-th
  root of unity. Then $\mathcal F$ is unitary and symmetric (it is the
  Discrete Fourier Transform matrix), so its inverse is its complex
  conjugate: in particular
  \[ (\wedge^r \mathcal F)^\vee \cdot (\wedge^l \mathcal F) = 
  \det F = (-1)^{r+l-1} I \]
  where $\bullet^\vee$ denotes transposition and $I$ is the identity matrix. 

  Hence $\wedge^r \mathcal F = (-1)^{r+l-1} \wedge^l \bar{\mathcal F}$,
  yielding:
  \[ \frac{Q_\lambda(u, u \in S)}{(r+l)^{r/2}} = 
  (-1)^{r+l-1} \eps(\lambda) \eps(S, T) 
  \frac{Q_{{\lambda^\ast}^T}(\bar{v}, v \in T)}{(r+l)^{l/2}} =
  (-1)^{r+l-1} \eps(\lambda) \eps(S, T) 
  \frac{Q_{\lambda^T}(v, v \in T)}{(r+l)^{l/2}}, \]
  where the second equality arises from the identity $\bar v^k = v^{r+l-k}$
  which holds for any $(r+l)$-th root $v$ of unity, and $\eps(S,T)$ denotes
  the signature of the permutation mapping $S$ to the $r$ first integers
  and $T$ to the next $l$ integers, and $\eps(\lambda)$ is the analogous
  permutation mapping $(r-1, \dots, 0) + \lambda$ to $(r-1, \dots, 0)$.
  The boxes of $\lambda$ describe the needed transpositions, so 
  $\eps(\lambda) = (-1)^{\abs\lambda}$.

  The identity we want follows, dividing by the case $\lambda=0$ which
  gives an identity for the Vandermonde determinant. 
\end{proof}

\begin{cor}
  The spaces $H^0(\SUpar{r}, \Lpar{l})$ and $H^0(\Utpar{l}, \Ltpar{r})$ 
  have equal dimensions for any $\vect\lambda$ whose size is a multiple of
  $rl$. 
\end{cor}

\begin{proof}
  We need to check that
  \[ \prod_{s \in S} (\zeta^s)^{-\abs{\vect\lambda}/r} 
  \prod_\lambda S_\lambda(\zeta^s, s\in S) =
  \prod_{t \in T} (\zeta^t)^{-\abs{\vect\lambda}/l} 
  \prod_{\lambda^T} S_\lambda(\zeta^t, s\in S) \] 
  which reduces to
  \[ \prod_{s \in S} (\zeta^s)^{-\abs{\vect\lambda}/r} = \prod_{t \in T}
  (\zeta^t)^{-\abs{\vect\lambda}/l} (-1)^{\abs{\vect\lambda}}. \]
  Now $\prod_{s \in S} (\zeta^s) \prod_{t \in T} (\zeta^t) = (-1)^{r+l-1}$
  and $-\abs{\vect\lambda}/l \equiv \abs{\vect\lambda}/r \mod (r+l)$.
  It now suffices to compare parities of $(r+l-1) \abs{\vect\lambda}/r$,
  which it the same as $(r+l-1) \abs{\vect\lambda}/l$ since their sum is
  $(r+l)(r+l-1) \frac{\abs{\vect\lambda}}{rl}$ with $\abs{\vect\lambda}$.
  Notice that $(r+l-1)rl/r = rl + l(l-1)$ has the same parity as $rl$,
  giving the result.
\end{proof}

The Witten correspondance between the Verlinde algebra and the quantum
cohomology of the Grassmannian will have a better geometric interpretation
if we use twisted Verlinde numbers, whose dimensions satisfy more
symmetries and gives more straightforward equalities.

\begin{prop}
  Let $\det: \Upar{r} \to \Jac(X)$ be the morphism associated to the
  $r$-th exterior power functor. Let $\Theta$ be a theta divisor on the
  Jacobian of $X$. It will also denote the associated line bundle and
  $\bar\Theta$ its pull-back by the involution $x \mapsto -x$. 

  Then if $\vect\lambda$ is a multiple of $rl$, 
  \[ \begin{aligned}
    h^0(\Upar{r}, &\Lpar{l} \otimes {\det}^\star \Theta) \\
    = &\sum_{S \in \mathfrak P_r(\Z/(r+l)\Z)} 
    \prod_{s \in S, t \notin S} \abs{2 \sin \pi\frac{s-t}{r+l}}^{g-1}
    \prod_{s \in S} (\zeta^s)^{-\abs{\vect\lambda}/r} 
    \prod_\lambda S_\lambda(\zeta^s, s\in S)
  \end{aligned} \]
\end{prop}

\begin{proof}
  Consider
  \[ \tau: \SUpar{r} \times \Jac(X) \to \Upar{r} \]
  again, and remember that $\tau^\star \Lpar{l} = \Lpar{l} \boxtimes
  \Theta^{rl}$. We also have $\tau^\star({\det}^\star \Theta) = \O
  \boxtimes r^\star \Theta$ which is numerically equivalent to $r^2
  \Theta$, thus they are interchangeable when computing with
  Grothendieck-Riemann-Roch theorem for example. This allows to write
  \[ \begin{aligned}
    h^0(\tau^\star (\Lpar{l} \otimes {\det}^\star \Theta) 
    &= r^{2g} h^0(\Upar{r}, \Lpar{l} \otimes {\det}^\star \Theta) \\
    &= h^0(\SUpar{r}, \Lpar{l}) h^0(\Jac X,\Theta^{rl}\otimes r^\star\Theta) \\
    &= h^0(\SUpar{r}, \Lpar{l}) h^0(\Jac X, \Theta^{rl+r^2}) \\
    &= r^g (r+l)^g h^0(\SUpar{r}, \Lpar{l}). \\
  \end{aligned} \]
  
  We can recover the desired formula by using \ref{Verlinde-formula}.
\end{proof}

\begin{cor}
  The spaces $H^0(\Upar{r}, \Lpar{l} \otimes {\det}^\star \Theta)$ and
  $H^0(\Utpar{l}, \Ltpar{r} \otimes {\det}^\star \bar\Theta)$ have equal
  dimensions. 
\end{cor}

\section{Symmetries associated to outer automorphisms}
\label{section:symmetries}

The affine Lie algebra of type $A_{r-1}^{(1)}$, which is
$\widehat{\mathfrak{sl}_r}$ has a Dynkin diagram with the shape of a
$r$-cycle: its automorphism group is the dihedral group $\Dcyr_r$. The
theory of affine Lie algebras \cite{Kac} identifies this group with the
outer automorphism group of the corresponding Kac-Moody algebra, and it is
known \citelist{\cite{FuchsSchlSchw} \cite{FuchsSchw}} that this group
induces correspondances between representations with different highest
weights and isomorphisms between the associated bundles of conformal blocks
in a very general setting.

In the particular $A_n^{(1)}$ case, these isomorphisms have simple
meanings. We will study the case of duality and rotations, which generate
the dihedral group, and show that strange duality morphisms are compatible
with these isomorphisms. 

\subsection{Conjugate moduli spaces}

\def\UBpar#1{\mathcal U_{X,\vect p}(#1, B)}
\def\SUBpar{\mathcal S\UBpar}
\def\Lpar#1{\mathcal L^{#1, \vect\lambda}}

The outer automorphisms induce isomorphisms between conformal blocks
which can be explained by isomorphisms between moduli spaces. Let
$\UBpar{r,d}$ be the moduli space of vector bundles with a parabolic
structure at the marked points given by full flags in the fibres, and
$\SUBpar{r,\Lcyr}$ be the analogous moduli space of bundles with fixed
determinant $\Lcyr$.

\begin{prop}
  For any $\vect\lambda$, there is canonical forgetful morphism
  $\UBpar{r,d} \to \Upar{r,d}$. This morphism is projective, so the
  conformal blocks can be expressed as sections of line bundles on
  $\Upar{r,d}$ as well as sections on $\UBpar{r,d}$ of their
  pull-backs.
\end{prop}

\begin{prop}[Duality]
  The duality automorphism, acting as a reflection on the regular Dynkin
  diagram induces isomorphisms:
  \[ V: \UBpar{r,d} \to \UBpar{r,-d} \]
  \[ V: \SUBpar{r,\Dcyr} \to \SUBpar{r,-\Dcyr} \]
  obtained by composition of the associated moduli functors with the
  duality functor.
\end{prop}

\begin{proof}
  The functor defining these automorphisms associates to a vector bundle
  its dual bundle with the orthogonal filtrations as seen in the previous
  sections. This functor is its own inverse, and identifies components
  of the moduli functor corresponding to dual determinants. 
\end{proof}

\begin{prop}[Rotations]
  The group $\Z/r\Z$ acts on the affine Dynkin diagram by $\rho: \varpi_i
  \to \varpi_{i+1}$ (with indices taken modulo $r$) and for any $\sigma \in
  (\Z/r\Z)^D$, there are isomorphisms
  \[ r_\sigma: \UBpar{r,d} \to \UBpar{r,d-\abs\sigma} \]
  \[ r_\sigma: \SUBpar{r,\Lcyr}\to\SUBpar{r,\Lcyr(-\sigma\cdot\vect p)} \]
  where $\abs{\sigma}$ is the length of $\sigma$ with respect to the
  generators of the form $\rho_p$ (which acts by $\rho$ at the point
  $p$), and $\sigma\cdot\vect p = \sum_p \sigma_p p$ (we identify
  $\SUBpar{r,\Lcyr}$ and $\SUBpar{r,\Lcyr(-rp)}$ using the tensor
  product by $\O(-p)$).
\end{prop}

\begin{proof}
  The functor $r_{\rho,p}$ representing the action of $\rho_p$ associates
  to a vector bundle $E$ with full flags
  \[ 0 \subset E_{p,1} \subset \cdots \subset E_{p,r} = E_p \]
  the vector bundle $F$ defined by the exact sequence 
  \[ 0 \to F \to E \to E_p/E_{p,r-1} \]
  There is an exact sequence 
  \[ 0 \to E_p/E_{p,r-1} = \Tor_1(\O_p, E_p/E_{p,r-1}) \to F_p \to E_p \to
  E_p/E_{p,r-1} \]
  so $F_p$ is an extension $i: F_p \to E_{p,r-1}$ by a line and has a
  canonical full flag 
  \[ 0 \subset E_p/E_{p,r-1} = i^{-1}(0) \subset \cdots \subset
  i^{-1}(E_{p,r-1}) = F_p \]

  To construct an inverse functor, we need to associate to a vector bundle
  $F$ with full flags an extension $E$ of $F$ by $\O_p$ given by an exact
  sequence 
  \[ 0 \to F \to E \to \O_p \to 0. \] 
  Such sequences are classified by a vector $f \in \Ext^1(\O_p, F)
  \simeq F_p$ which describes the behaviour of the fibre at $p$:
  there is an exact sequence
  \[ 0 \to \O_p \xrightarrow{f} F_p \to E_p \to \O_p \to 0. \]
  
  The locally free extensions are then classified by $\mathbb P(F_p)$
  but the only one yielding compatible flag structures is the point
  corresponding to the line $F_{p,1}$.

  These constructions extend to families of vector bundles and thus give
  isomorphisms between the moduli stacks. Iterating these functors gives
  the action of the full automorphism group.
\end{proof}

\begin{rem}
  This construction is related to Simpson's description (see
  \cite{Simpson}) of parabolic vector bundles as vector bundles with a
  \emph{periodic} filtration
  \[ \cdots \to E_{r-1}(-p) \to E(-p) \to E_1 \to \cdots \to E_{r-1} \to E
  \to E_1(p) \to \cdots E_{r-1}(p) \to E(p) \to \cdots \]

  The action of $\rho_p$ corresponding to shifts of this filtration,
  so $\rho_p$ will also be called the shift automorphism.
\end{rem}

These isomorphisms descend to isomorphisms between moduli spaces of vector
bundles with parabolic structures of type $\vect\lambda$ where
$\vect\lambda$ is a labelling of the marked points by weights of
$\mathfrak{sl_r}$. The group $(\Dcyr_r)^D$ acts by $\vect\sigma(\vect\lambda) =
(\sigma_p(\lambda_p))_{p \in D}$.

\begin{prop}
  There are isomorphisms
  \[ \Upar{r,d} \to \mathcal{U}_{\Xp}(r, -d, \vect\lambda^\ast) \]
  and 
  \[ \SUpar{r,\Dcyr} \to \mathcal{SU}_{\Xp}(r, -\Dcyr,
  \vect\lambda^\ast) \] 
  making the following diagrams commute
  \[ \xymatrix{
    \UBpar{r,d} \ar[r] \ar[d] & \UBpar{r,-d} \ar[d] \\
    \Upar{r,d} \ar[r] & \mathcal{U}_{\Xp}(r, -d, \vect\lambda^\ast) \\
  } \quad \xymatrix{
    \SUBpar{r,\Dcyr} \ar[r] \ar[d] & \SUBpar{r,-\Dcyr} \ar[d] \\
    \SUpar{r,\Dcyr} \ar[r] & \mathcal{SU}_{\Xp}(r,-\Dcyr, \vect\lambda^\ast) \\
  } \]
\end{prop}

\begin{prop}
  Let $0 < k < r$ be the smallest integer such that
  $\pair{\lambda_p,\alpha_{r-k}} > 0$ (where $\alpha_{r-k}$ are the simple
  roots of $\mathfrak{sl}_r$), there are morphisms
  \[ r_p: \Upar{r,d} \to \mathcal{U}_{\Xp}(r, d-k, \rho_p^k(\vect\lambda)) \]
  and 
  \[ r_p: \SUpar{r,\Dcyr} \to \mathcal{SU}_{\Xp}(r, \Dcyr(-kp),
  \rho_p^k(\vect\lambda)) \]
  making the following diagrams commute
  \[ \xymatrix{
    \UBpar{r,d} \ar[r] \ar[d] & \UBpar{r,d-1} \ar[d] \\
    \Upar{r,d} \ar[r] & \mathcal{U}_{\Xp}(r, d-1, \rho_p^k(\vect\lambda)) \\
  } \quad \xymatrix{
    \SUBpar{r,\Dcyr} \ar[r] \ar[d] & \SUBpar{r,\Dcyr(-p)} \ar[d] \\
    \SUpar{r,\Dcyr} \ar[r] & \mathcal{SU}_{\Xp}
    (r,\Dcyr(-p),\rho_p^k(\vect\lambda)) \\
  } \]
  as well as preserving the inclusions of moduli spaces of bundles with
  fixed determinant in moduli spaces of bundles with varying determinant.
\end{prop}

\subsection{Isomorphisms between spaces of conformal blocks}

\begin{prop}
  Consider the duality isomorphism:
  \[ V: \Upar{r,d} \to \mathcal{U}_{\Xp}(r, -d, \vect\lambda^\ast) \]
  The pull-back of line-bundles is given by the following formula:
  \[ V^\star \mathcal L^{l,\vect\lambda^\ast}_{M^\vee \otimes K_X}
  \simeq \Lpar{l}_M \]
  As the consequence, there is a canonical isomorphism between 
  $H^0(\mathcal{U}_{\Xp}(r, -d, \vect\lambda^\ast),
  \mathcal L^{l,\vect\lambda^\ast}_{M^\vee \otimes K_X})$ and 
  $H^0(\Upar{r,d}, \Lpar{l}_M)$ for any rank $l$ vector bundle $M$ with a
  parabolic structure of type $\vect\lambda^T$ and parabolic degree
  $d/r+1-g$. 
\end{prop}

\begin{proof}
  The identity follows from Serre duality: if $E$ is a family of
  parabolic vector bundles over $X$ parameterised by $S$, the direct
  image by $\pi: X \times S \to S$ satisfies
  \[ R\pi_\star (E \otimes M)_{\rm par} = R\pi_\star 
    (K_X \otimes E^\vee \otimes M^\vee)_{\rm par} [1]^\vee \]
  (the extension of Serre duality to parabolic bundles follows from
  the properties of dual parabolic structures).
\end{proof}

\begin{prop}
  We suppose $\lambda^{(r)} = 0$, so $\sigma(\lambda)$ can be defined
  at the level of Young diagrams. The shift morphism associated to an
  elementary rotation $\sigma$ at a point $p$
  \[ r_\sigma: \UBpar{r,d} \to \UBpar{r,d-1} \]
  transform line bundles according to the formula
  \[ r_\sigma^\star \mathcal L_M^{l,\sigma(\vect\lambda)} \simeq
  \Lpar{l}_M \]
  and induce isomorphisms between $H^0(\Upar{r,d},\Lpar{l}_M)$ and 
  $H^0(\mathcal{U}_{\Xp}(r, d-1, \mathcal L^{l,\sigma(\vect\lambda)}_M)$
  for any rank $l$ vector bundle $M$ with its parabolic structure of type
  $\vect\lambda^\ast \equiv \sigma(\vect\lambda)^\ast$.
\end{prop} 

\begin{proof}
  If $E$ is a family of vector bundles over $X$, and $E^\sigma$ is the
  vector bundle obtained by the functor associated to $\sigma$, the
  isomorphism between the determinant line bundles is induced by the
  canonical isomorphism between $(E \otimes M)_{\rm par}$ and
  $(E^\sigma \otimes M_\sigma)_{\rm par}$ where $M_\sigma$ is the
  vector bundle $M$ equipped with its parabolic structure of type
  $\lambda + \varpi_l$ (with is equivalent to the structure of type
  $\lambda$).

  We now observe that elements $\Hom_{\rm par}(M^\vee, E)$ must
  actually be maps to $E^\sigma$ (because $\lambda^{(r)} = 0$), and
  can be identified with elements of $\Hom_{\rm par}(M^\vee_\sigma,
  E^\sigma)$.
\end{proof}

In the cases where the weights used have non trivial stabilisers, there can
be automorphisms of the space of conformal blocks. 

\begin{prop}
  Let $\vect\lambda$ be a system of weights of $\mathfrak{sl}_r$ and let
  $\phi \in (\Z/r\Z)^D$ be a choice of a rotation at each point
  such that $\phi_p(\lambda_p) = \lambda_p$ and $\sum \phi_p = 0$. 
 
  Then $\phi$ induces an isomorphism $\UBpar{r,d} \simeq \UBpar{r,d+kr}$ 
  and the choice of a $r$-th root of $\O(\sigma\cdot\vect p)$ turns $\phi$ into
  an automorphism of $\UBpar{r,d}$ commuting with the determinant morphism 
  $\UBpar{r,d} \to \Pic^d X$. 

  In particular, $\phi$ defines automorphisms of the spaces
  $H^0(\SUpar{r,\Dcyr}, \Lpar{l})$ for $\Dcyr \in \Pic^d X$. 
\end{prop}

\subsection{Behaviour of rank-level duality}

\begin{prop}
  The morphisms between spaces of conformal blocks defined by the duality
  functor commute with the strange duality morphisms defined in section 1.
\end{prop}

\begin{proof}
  This is a consequence of the compatibility of tensor products with
  duality. 
\end{proof}

To study the behaviour of rank-level duality with respect to the outer
automorphisms, we need more subtle operations.

\begin{defn}
  Let $\lambda$ be a Young diagram in a $r \times l$ rectangle. The
  \emph{string} of $\lambda = (a_r, \dots, a_1)$ (where $a_i$ is
  non-increasing) is the sequence $(x_1, \dots, x_{r+l})$ where $x_i = R$
  if $i = a_k+k$ for some $k$ and $x_i = L$ otherwise (here $R$ and $L$ are
  abstract symbols). 
\end{defn}

\begin{defn}[String parabolic structure]
  Let $E$ and $F$ be vector spaces with full flags and let $\lambda$ be a
  Young diagram inside a $r \times l$ rectangle. Denote by 
  $(E_i)_{i=1}^r$ and $(F_j)_{j=1}^l$ the increasing filtrations on $E$ and
  $F$. The space $\Hom(F,E)$ has a canonical
  parabolic structure of level $r+l$ defined by the subspaces 
  \[ G_k = \Hom(F/F_{k_L}, E_{k_R}) \subset \Hom(F,E) \]  
  for $k \in \set{1, \dots, r+l}$ where $k_R$ is the number of $i > k$ such
  that $x_i = R$ and $k_L$ is the number of $j \leq k$ such that $x_j = L$.
\end{defn}

\begin{rem}
  Notice that $E_{k_R}$ decreases and that $F_{k_L}$ increases, so that
  $G_k$ is a decreasing filtration. 
\end{rem}

Given two parabolic bundles of types $\vect\lambda$ and $\vect\lambda^T$,
with a choice of compatible full flags at each point, we now endow $E
\otimes F$ with the structure $\abs{\vect\lambda}^+$ which consists in the
natural structure $\abs{\vect\lambda}$ with the additional data of $E_{p,1}
\otimes F_p$ at $p$, where $E_1$ is the chosen hyperplane at $p$.

\begin{prop}
  Let $\sigma$ be the elementary rotation at the point $p$, and suppose 
  $\lambda_p^{(r)} = 0$. Then $\abs{\sigma(\lambda_p)} = \abs{\lambda_p} +
  l$. We have natural morphisms
  \[ \xymatrix@C-=3ex{
    \mathcal U(r, d, B) \ar[d]^{r_\sigma} \ar@{}[r]|{\times}
    & \mathcal U(l, \dcyr, B) \ar@{=}[d] \ar[r]^(.3){\tau_{\lambda}}
    & \mathcal U(rl, ld+r\dcyr, \abs{\vect\lambda}^+) 
    \times \Pic^{\dcyr-d}(X) \ar[d]^{r_{\sigma^l}}\\
    \mathcal U(r, d-1, B) \ar@{}[r]|(.55){\times} & \mathcal U(l, \dcyr, B)
    \ar[r]^(.25){\tau_{\sigma(\lambda)}} & \mathcal U(rl, l(d-1)+r\dcyr,
    \abs{\sigma{\vect\lambda}}) \times \Pic^{\dcyr-d+1}(X) \\
  } \]
  Here $\tau_\lambda$ is the tensor product with respect to type
  $\lambda$-structures with additional $+$-structure. We also have the
  corresponding pull-back functors between line bundles
  \[ \xymatrix@C-=2ex{
    \mathcal L_{F_0 \otimes L}^{l,\vect\lambda} \otimes
    {\det}^\star (\det F_0 - \Theta) \ar@{}[r]|{\boxtimes}
    & \mathcal L_{E_0 \otimes L}^{r,\vect\lambda^T} \otimes
    {\det}^\star (\det E_0 + \Theta) 
    & \mathcal L_L^{1,\abs{\vect\lambda}} 
    \boxtimes \Theta \ar[l] \\
    \mathcal L_{F_0 \otimes L}^{l,\sigma(\vect\lambda)} \otimes
    {\det}^\star (\det F_0 - (\Theta + p)) \ar[u] \ar@{}[r]|{\boxtimes}
    & \mathcal L_{E_0' \otimes L}^{r,\vect{\lambda^T}} \otimes
    {\det}^\star (\det E_0' + (\Theta+p)) \ar@{=}[u]
    & \mathcal L_L^{1,\abs{\sigma(\vect{\lambda})}} 
    \boxtimes (\Theta + p) \ar[l] \ar[u]^{\tau^\star} \\
  } \]
  where $E_0 \in \mathcal U(r, d, \vect\lambda)$ and $F_0 \in\mathcal
  U(l, \dcyr, \vect{\lambda^T})$ and the corresponding pull-back
  isomorphisms between spaces of sections commute. Moreover, these
  isomorphisms preserve the canonical sections of the parabolic determinant
  bundles, so they commute with the strange duality morphisms defined in
  section \ref{subsection:arbitrary-degree}.
\end{prop}

\begin{proof}
  The commutativity results from the fact that these line bundles are all
  determinant bundles of the same vector bundle (up to pull-back) of rank
  $rl$ and slope $g-1$.
\end{proof}

\section{Associated enumerative problems}

\subsection{Enumerative interpretation of the strange duality}
\label{subsection:enumerative-interpretation}

The exceptional symmetries shown in the previous sections allow to prove
strange duality for special degrees or weights without loss of generality.

\begin{prop}
  Let $d$ and $\dcyr$ be integers, and $\vect\lambda$ be a system of 
  possibly empty Young diagrams such that 
  \[ \abs{\vect\lambda} +ld+r\dcyr  \equiv 0 \mod rl. \]

  Up to adding enough points wearing empty Young diagrams, there
  exists a sequence of rotations of Young diagrams inducing maps
  $\UBpar{r,d} \to \UBpar{r,d'}$ and $\UBpar{l,\dcyr} \to
  \UBpar{l,\dcyr'}$ such that the resulting system of Young diagrams
  $\vect\nu$ and the corresponding degrees $d'$, $\dcyr'$ and the new
  number of marked points $n'$ satisfy:
  \begin{itemize}
  \item $d'+rn'$ is an arbitrary large multiple of $r$,
  \item $\dcyr' = 0$,
  \item $\abs{\vect\nu^\ast} = n'rl + ld'+r\dcyr' + rl(1-g)$.
  \end{itemize}
\end{prop}

\begin{proof}
  Since rotating an empty Young diagram gives $l\varpi_1$,
  $\abs{\vect\lambda}$ increases by $l$ so adding enough points with empty
  diagrams and rotating them can make the total size of the diagrams
  arbitrarily large. Each step of this process reduces $d$ by $1$, and
  preserves $\dcyr$, so we can obtain an value $d'$ which is multiple of
  $r$.

  Applying the previous argument to the other moduli space, we can
  increase or decrease $\dcyr$ to zero while not changing $d$ (see the
  previous section for the behaviour of degrees under rotation of the
  diagrams), using rotations or their inverse of empty diagrams.

  The third equality automatically holds modulo $rl$. Tensoring by a line
  bundle of degree $s$ gives an isomorphism between $\Upar{r,d'}$ and
  $\Upar{r,d'+rs}$ increasing $ld'$ by $rl$: choosing an appropriate $s$ we
  obtain an equality for the third relation (this can be achieved by
  applying full rotations to empty Young diagrams). Since $\abs{\vect\nu}$
  can be arbitrarily large, $d'+rn'$ can be as large as needed.
\end{proof}

Recall there is a tensor product morphism
\[ \widetilde\tau: \Upar{r,d} \times \Utpar{l,\dcyr} \to \Uppar{rl,
  ld+r\dcyr} \times \Pic^{\dcyr-d} X \] 
which pulls the canonical section of $\Lppar{1}_\Lcyr \boxtimes \Theta$ to
a strange duality pairing $\widetilde{SD}$ which is a section of
$(\Lpar{l}_{F_0 \otimes \Lcyr} \otimes {\det}^\star (\det F_0 - \Theta))
\boxtimes (\Ltpar{r}_{E_0 \otimes \Lcyr} \otimes {\det}^\star (\det E_0 +
\Theta))$, where $\Lcyr$ is a line bundle on $X$ of degree
  \[ g-1 - \frac{(ld+r\dcyr) + \abs{\vect\lambda}}{rl}. \]

The previous proposition allows to replace the moduli spaces by moduli of
bundles with different degrees. The rank-level duality between the
conformal blocks is then equivalent to the rank-level duality in the new
setting. 

\begin{prop}
  One can replace $d$, $\dcyr$ and $\vect\lambda$ by values $d'$,
  $\dcyr'=0$, $\vect\nu$ satisfying the conditions above and define
  canonical isomorphisms between conformal blocks, pulling the rank-level
  duality morphism to the new one.

  In particular, proving rank-level duality in one of these contexts
  is equivalent to proving it in the other one.
\end{prop}

\begin{rem}
  After this operation, the line bundle used to define the duality has
  degree zero, so without loss of generality we can suppose it is trivial.
\end{rem}

\begin{cor}
  The rank-level symmetry between Verlinde numbers holds in arbitrary
  degree in the case 
  \[ \abs{\vect\lambda} + ld+r\dcyr \equiv 0 \mod rl. \]
\end{cor}

We now suppose the conditions above are satisfied by $\vect\lambda$, and
replace $d$ par a very large number so that
$\abs{\vect\lambda^\ast} = ld + rl(n+1-g)$ is a multiple of $rl$ and call $q$ the
common dimension of the vector spaces of conformal blocks.

\begin{lemma}
  Suppose there exist $q$ pairs of vector bundles $(A_i, B_i)$ of degree
  $(d,0)$ and ranks $r$ and $l$, with parabolic structure of type
  $\vect\lambda$ and $\vect\lambda^T$ at points $\vect p$, and a line
  bundle $M$ of degree $g-1+d$ defining a theta divisor $\Theta_M$ on
  $\Pic^{-d}(X)$ satisfying the following properties:
  \begin{itemize}
  \item $H^0((A_i \otimes B_j)_{\rm par}) = 0$ iff $i = j$,
  \item $H^0(\det A_i^\vee \otimes \det B_j \otimes M) = 0$.
  \end{itemize}

  Then the strange duality morphism $\widetilde{SD}$ induced by the
  canonical section of 
  \[ H^0(\Uppar{rl}, \Lppar{1}) \otimes H^0(\Pic^{-d} X, \Theta_M) \]
  is an isomorphism between 
  \[ H^0(\Upar{r,d}, \Lpar l_{F_0} \otimes 
  {\det}^\star(\det F_0-\Theta_M)) \]
  and
  \[ H^0(\Utpar{l}, \Ltpar r_{E_0} \otimes 
  {\det}^\star(\det E_0+\Theta_M)). \]
\end{lemma}

\begin{rem}
  Note the second condition can easily be fulfilled by choosing $M$ in the
  complement of a finite number of theta divisors (defined by the
  non-satisfaction of this condition).
\end{rem}

\begin{proof}
  From the definition of the canonical section of the determinant bundle,
  the properties of the pairs amount to the non-vanishing of the strange
  duality pairing at points $(A_i, B_j)$ if and only if $i=j$, which means
  $(A_i)$ and $(B_j)$ are linearly independant as elements of the dual of
  the spaces of generalised theta functions, and in duality for the strange
  morphism. Since there are $q$ couples of such bundles, they are actually
  bases and the strange morphism is a perfect pairing.
\end{proof}

\subsection{Geometric realisations of the Verlinde numbers}

The Vafa-Intriligator formula, which is formulated and proved e.g. in
\cite{MarianOprea-virtual}, expresses intersection numbers on $\Quot$
schemes (parameterising quotients of $\O_X^{r+l}$) and the structure
coefficients of the quantum cohomology of the Grassmannian. As in the
non-parabolic case, it yields values very similar to the Verlinde numbers.
For example, the factorisation formulae for the Verlinde algebra (see
\cite{Beauville-verlinde}) correspond to the factorisation properties of
Gromov-Witten invariants as stated in \cite{SiebertTian}.

This is a geometric interpretation of the isomorphism between the Verlinde
algebra of $\mathfrak{gl}_r$ at level $l$ and the quantum cohomology of
$\Gr(r,r+l) \simeq \Gr(l,r+l)$. It will in the end give a construction of
points of the moduli spaces representing the strange duality.

Let $r$, $l$, $d$, be integers, and $\vect\lambda$ be a system of Young
diagrams such that $\abs{\vect\lambda} = ld + rl(n+1-g)$. We assume that $d$
and $n$ are multiples of $2rl(r+l)$ (up to adding more points), according
to the previous section. We set $\vect\mu = \vect\lambda^{\ast}$.

Let $V = L^{\oplus r+l}$ be a vector bundle on $X$ with rank $r+l$ and
degree $rn+d$. Let $Q = \Quot_{r,d+rn}(V)$ be the fine moduli scheme
parameterising quotients of degree $d+rn$ and rank $r$ of $V$. Universality
yields on $Q \times X$ an exact sequence of coherent sheaves
\[ 0 \to \mathcal E^\vee \to V \to \mathcal F(D) \to 0 \]
where $\mathcal E$ is a flat family over $Q$ of rank $l$ and degree $0$ 
vector bundles on $X$, and $\mathcal F$ is a flat family of rank $r$ and
degree $d$ coherent sheaves on $X$. 

The tensor product by $L^\vee$ induces an isomorphism of $Q$ with
$\Quot_{r,\frac{nrl+ld}{r+l}}(\O^{r+l})$ so the known results on this Quot
scheme also apply to $Q$. In particular, $Q$ is a compactification of 
\[ \Mor_{\frac{nrl+ld}{r+l}}(X, \Gr(l, r+l)), \] 
\noindent which is identified to the locus of $Q$ representing locally free
quotients of $V$. 

If $rn+d$ is large enough, it is known that $Q$ is a reduced irreducible
scheme \cite{BDW} of dimension $ld+rl(n+1-g) = \abs{\vect\mu}$.

Let $\iota: Q \to Q \times X$ be the constant section mapping $q \in Q$ to
the point $(q,p)$ where $p$ is a fixed point. Let $a_k$ (for $k$ between
$1$ and $l$) be the Chern classes of $\iota^\star (\mathcal E \otimes
L^\vee)$. Then the numerical equivalence class of $a_k$ does not depend on $p$
for a generic choice of this point.

Let $\lambda$ be a Young diagram with $\lambda^T \leq (r^l)$.  The
characteristic class $a_\lambda$ is defined as $R_{\lambda^T}(1, a_1,
\dots, a_l)$ where $R_{\lambda^T}$ is the polynomial computing the Schur
polynomial in terms of elementary symmetric functions.

\begin{prop}[see \cite{Fulton-schubert}]
  The characteristic class $a_{\mu_p}$ is represented by the subvariety
  $Z_p$ of $Q$ parameterising morphisms from $X$ to $\Gr(l, r+l)$ mapping
  $p$ into a given Schubert subvariety $Y_p$ of type $Y_{\mu_p^T}$, which has
  codimension $\abs{\mu_p}$.
\end{prop}

\begin{prop}[see \cite{Bertram-schubert}]
  There exist a choice of Schubert varieties such that the $Z_p$ intersect
  properly, on the smooth locus of $Q$ parameterising stable locally free
  quotients of $V$. 
\end{prop}

\begin{prop}
  If $a_{\vect\mu} = \prod_p a_{\mu_p}$, we have for any parabolic vector
  bundle $F_0$ of rank $r$ and degree $d$,
  \[ h^0(\Utpar{l,0}, \Ltpar{r}_{F_0} \otimes {\det}^\star \bar\Theta) =
  \int_Q a_{\vect\mu}. \]
\end{prop}

\begin{proof}
  The intersection number above can be computed by the Vafa-Intriligator
  formula, as stated in \cite{MarianOprea-virtual}. We obtain: 
  \[ \int_Q a_{\vect\mu} = (r+l)^{l(g-1)}
  \sum_{T \in \mathfrak P_l(\Z/(l+r)\Z)}
    \prod_p S_{\lambda_p}(\zeta^T) \pa{ \prod \zeta^t 
    \prod_{t \neq u \in T} (\zeta^t - \zeta^u)}^{1-g} \]
  where $T$ goes through order $l$ subsets of $\Z/(r+l)\Z$. We have the
  additional property 
  \[ \prod_{t \neq u \in T} (\zeta^t - \zeta^u) = 
  \prod_{t \neq u \in T} \abs{\zeta^t - \zeta^u} 
  (\prod_{t \in T} \zeta^t)^{l-1} \]
  since $(\zeta^s-\zeta^t)(\zeta^t-\zeta^s) = \zeta^s\zeta^t
  \abs{\zeta^s-\zeta^t}^2$. 
  Factoring out $(\prod_{t \in T} \zeta^t)$ at the power $l(1-g)$ (which is
  congruent to $\abs{\vect\mu}/r \equiv -\abs{\vect\mu}/l$ modulo $(r+l)$),
  there remains:
  \[ \int_Q a_{\vect\mu} =
  (r+l)^{l(g-1)} \sum_T \prod_p S_{\mu_p}(\zeta^T) 
  (\prod_{t \in T} \zeta^t)^{-\abs{\vect\mu}/l}
  \prod_{t \neq u \in T} \abs{\zeta^t - \zeta^u}^{1-g} \]
  which equals $h^0(\Utpar{l}, \Ltpar{r})$ by the Verlinde formula
  \[ h^0(\Utpar{l}, \Ltpar{r} \otimes {\det}^\star \bar\Theta) =
  (r+l)^{l(g-1)} \sum_T (\prod \zeta^T)^{-\abs{\vect\lambda}/r}
  \prod_\lambda S_\lambda(\zeta^T)
  \prod_{t \neq u \in T} \abs{\zeta^t - \zeta^u}^{1-g} \]
\end{proof}

The intersection cycle $a_{\vect\mu}$ can be represented by a reduced
zero-dimensional subscheme $Z$ that counts the number of points of $Q$
representing morphisms from $X$ to $\Gr(l,r+l)$ mapping the marked points
$p$ of $X$ to the Schubert cells $Y_p$. But we saw that a point of a
Schubert cell $Y_{\mu^T_p}$ inherits a natural parabolic structure of type
$\mu$.

A point of $Z$ is an exact sequence
\[ 0 \to E^\vee \to V \to F(D) \to 0 \]
such that $E_p^\vee$ is a $l$-dimensional subspace of $V_p$ chosen in a
particular Schubert variety $Y_p$: this endows $E^\vee$ with a type
$\vect\mu^T$ parabolic structure. By duality, $E$ has rank $l$ and degree
$0$, and is equipped with a type $\vect\lambda^T$ parabolic structure.

We also showed that $F(D)^\vee = (E^\vee)^\perp \subset V^\vee$ inherited a
natural $\vect\mu$ structure, giving $F$ a type $\vect\lambda$
parabolic structure.

\begin{prop}
  The points of $Z$ parameterise $q$ pairs of vector bundles $(E_i, F_i)$
  of degree $(-d,0)$ and ranks $r$ and $l$, with parabolic structure of type
  $\vect\lambda$ and $\vect\lambda^T$ at points $\vect p$, such that 
  ${H^0(E_i \otimes F_j(D))_{\rm par} = 0}$ iff $i = j$.
\end{prop}

\begin{proof}
  Consider the exact sequences of parabolic bundles
  \[ 0 \to \mathcal E_i^\vee \to V \to \mathcal F_i(D) \to 0 \]
  given by the intersection points of Schubert cells in the previous
  section. They are the points of a reduced zero-dimensional scheme
  (defined by the inverse image of Schubert varieties in the $\Quot$
  scheme). Its tangent space can be written as $\Hom_{\rm par}(\mathcal
  E_i^\vee, \mathcal F_i(D))$, using the previous computations of the tangent
  space of Schubert cells, and the fact that $Z$ was a proper
  intersection. This tangent space vanishes since $Z$ is smooth of
  dimension zero.

  Moreover, replacing $F_i$ by $F_j$ breaks the exactness of the sequence,
  and the composite $E_i^\vee \to V \to F_j$ gives a nonzero morphism (if
  it were not the case, there would be a morphism from $E_i^\vee$ to
  $E_j^\vee$, and it would be injective since $E_i^\vee \to V$ is
  injective, and since $E_i$ and $E_j$ have the same degree this is
  impossible). 
\end{proof}

This solves the enumerative problem we set up to prove strange duality in
section \ref{subsection:enumerative-interpretation}. 

\begin{thm}[Rank-level duality of $GL_n$ conformal blocks]
  Let $(\Xp)$ be a smooth projective curve with marked points, and
  $\vect\lambda$ be a system of Young diagrams.

  Let $d$ and $\dcyr$ be integers such that
  \[ \abs{\vect\lambda} + ld+r\dcyr \]
  is divisible by $rl$.

  Consider the tensor product
  \[ \widetilde\tau: \Upar{r,d} \times \Utpar{l,\dcyr} \to
  \Uppar{rl, ld+r\dcyr} \times \Pic^{\dcyr-d} X \]
  mapping $(E,F)$ to their parabolic tensor product $(E \otimes F, \det
  E^\vee \otimes \det F)$. Let $\Lcyr$ be a line bundle with slope
  $g-1-(ld+r\dcyr+\abs{\vect\lambda})/(rl)$ and $\Theta$ be a generic
  Theta divisor on $\Pic^{\dcyr-d} X$. Then choosing a base point $(E_0,
  F_0)$ defines a divisor $\overline\Theta_d$ on $\Pic^d X$ pulling back by
  $x \to \det F_0 - x$ and $\Theta_\dcyr$ on $\Pic^\dcyr X$ pulling back by
  the action of $\det E_0^\vee$.

  Then
  \[ \widetilde\tau^\star \Lppar{1}_\Lcyr \boxtimes \Theta =
  (\Lpar{l}_{F_0 \otimes \Lcyr} \otimes {\det}^\star \overline\Theta_d)
  \boxtimes
  (\Ltpar{r}_{E_0 \otimes \Lcyr} \otimes {\det}^\star \Theta_\dcyr) \]

  Moreover there exists a section in $H^0(\Uppar{rl, ld+r\dcyr},
  \Lppar{1}_\Lcyr) \otimes {\det}^\star \Theta$, canonical up to a scalar,
  whose pull-back defines a perfect pairing
  \[ \widetilde{SD}: H^0(\Upar{r,d}, \Lpar{r}_{F_0 \otimes \Lcyr} 
  \otimes {\det}^\star \overline\Theta_d)^\vee \otimes
  H^0(\Utpar{l,\dcyr}, \Ltpar{l}_{E_0 \otimes \Lcyr} 
  \otimes {\det}^\star \Theta_\dcyr)^\vee \to \C. \]
\end{thm}

\begin{thm}[Rank-level duality of $SL_n/GL_n$ conformal blocks]
  With the same hypotheses, let $\Dcyr$ be a degree $\dcyr$ line bundle on
  $X$. The tensor product of parabolic bundles defines a morphism
  \[ \tau: \Upar{r,d} \times \SUtpar{l,\Dcyr} \to 
  \Uppar{rl,ld+r\dcyr}. \]   
  
  If $F_0$ is a given point in $\SUtpar{l,\Dcyr}$ and $\Lcyr$ a line bundle
  with appropriate slope as before, we have
  \[ \tau^\star \Lppar{1}_\Lcyr = \Lpar{l}_{F_0 \otimes \Lcyr} \boxtimes 
  \Ltpar{r} \]
  which does not depend on $F_0$ (since it always has determinant
  $\Dcyr$).

  The pull-back of the canonical section of $\Lppar{1}$ defines 
  a perfect pairing 
  \[ SD: H^0(\Upar{r,d}, \Lpar{l}_{F_0 \otimes \Lcyr})^\vee \otimes
  H^0(\SUtpar{l,\Dcyr}, \Ltpar{l})^\vee \to \C. \]
\end{thm}

\begin{rem}
  The rank-level duality for nonzero degree bundles proved in
  \cite{MarianOprea} corresponds to the case $\vect\lambda=0$
  and $rd+l\dcyr=0$.
\end{rem}

\section{Rank-level duality for Sp$_2$-bundles and Sp$_n$-bundles}

Here we consider moduli stacks of principal $G$-bundles with parabolic
structure (see \cite{LaszloSorger} for example). If $\lambda$ is a Young
diagram in a $r \times l$-rectangle, we define parabolic structures of type
$\lambda$ in a symplectic vector space of dimension $2r$ to be a partial
flag of isotropic subspaces whose dimensions are given by the row lengths
of $\lambda^T = (b_1, \dots, b_l)$. These structures correspond to choices
a parabolic subgroup conjugate to $P_\lambda$ in $Sp_{2r}$.

Let $\widetilde\lambda$ be the Young diagram in a rectangle of size $(2r
\times 2l)$ of the following form
\[ \begin{bmatrix}
  \lambda & 0 \\
  1 & \lambda^\ast \\
\end{bmatrix} \]

If $0 \subseteq E_1 \subseteq \cdots \subseteq E_l \subset E$ is an
isotropic partial flag in $E$, we define an associated flag 
\[ 0 \subseteq E_1 \subseteq \cdots \subseteq E_l \subseteq
E_l^\perp \subseteq \cdots \subseteq E_1^\perp \subseteq E \]
which determines a parabolic structure of type $\widetilde\lambda$ on 
$E$ as a plain vector space. 

As in the $GL_n$ case, we choose a system of Young diagrams $\vect\lambda$
labelling the marked points. Let $\mathcal M_\Xp(r, \vect\lambda)$ be the
moduli stack of symplectic vector bundles with parabolic structures of type
$\vect\lambda$. We also define an analogous moduli stack $\mathcal
M'_\Xp(r, \vect\lambda)$ for vector bundles equipped with a symplectic form
with values in $K_X(D)$ and a parabolic structure of type $\vect\lambda$.

If $\lambda$ is a Young diagram, we denote by $\overline{\lambda}$ the
half-integer weight of $\mathfrak{so}_{rl}$ which has the same parity as
$\abs{\lambda}$. Let $\mathcal Q'_\Xp(4rl,\overline{\vect\lambda})$ be the
moduli stack of oriented orthogonal bundles of rank $4rl$ whose quadratic
form takes values in $K_X(D)$, equipped with a parabolic structure of type
$\overline{\vect\lambda}$.

\begin{prop}[\cite{Abe}]
  The tensor product of symplectic bundles of ranks $2r$ and $2l$ with
  parabolic structures of type $\vect\lambda$ and $\vect\lambda^T$ is an
  orthogonal bundle with a canonical parabolic structure of type
  $\overline{\vect\lambda}$ coinciding with the usual parabolic structure
  for the tensor product of plain vector bundles.
\end{prop}

\begin{prop}
  We have a morphism of moduli stacks 
  \[ \mathcal M_\Xp(r, \vect\lambda) \times \mathcal M'_\Xp(l, \vect\lambda)
  \to \mathcal Q'_\Xp(rl, \overline{\vect\lambda}) \]
  and a commutative diagram :
  \[ \xymatrix{
    \mathcal U_\Xp(r, \widetilde{\vect\lambda}) \times 
    \mathcal U_\Xp(l, l(g-1-n/2), \widetilde{\vect\lambda^T}) \ar[d] & \\
    \mathcal M_\Xp(r, \vect\lambda) \times 
    \mathcal M'_\Xp(l, \vect\lambda^T) \ar[d] \ar[r]
    & \mathcal Q'_\Xp(rl, \overline{\vect\lambda}) \ar[d] \\
    \mathcal SU_\Xp(2r, \widetilde{\vect\lambda}) \times 
    \mathcal SU_\Xp(2l, \widetilde{\vect\lambda^T}) \ar[r] &
    \mathcal SU_\Xp(4rl, |\widetilde{\vect\lambda}|) \\
  } \]
  where morphisms from the first line to the second one map respectively
  vector bundles $E$ to symplectic bundles $E \otimes E^\vee$ and vector
  bundles $F$ to symplectic bundles $F(D) \otimes K_X \otimes F^\vee$.
\end{prop}

The moduli spaces of parabolic symplectic bundles carry natural line
bundles denoted by $\mathcal L^{l,\lambda}$ and $\mathcal
L^{r,\lambda^T}$. It is known (see \cite{OxburyWilson} for non-parabolic
bundles) that the spaces of sections of these line bundles have the same
dimensions. Moreover, there is a natural morphism from $\mathcal
Q'_\Xp(4rl,\overline{\vect\lambda})$ to $\mathcal Q'_X(4rl)$ : the moduli
stack of vector bundles of rank $4rl$ with a nondegenerate quadratic form
with values in $K_X$. This morphism represents the functor associating to a
parabolic orthogonal vector bundle the orthogonal vector bundle of sections
going through the chosen maximal isotropic subspaces (notice that the
target of the quadratic form changes). The Pfaffian line bundle on the
moduli stack $\mathcal Q'_X(4rl)$ has a canonical section, which is a
square root of the canonical section of the determinant bundle : it
vanishes on the locus of bundles with nonzero sections.

If the case where $l = 2$, the isomorphism between $Sp_2$ and $SL_2$ allows
to identify more vector spaces (see also \cite{MarianOprea-survey}) :
\begin{prop}
  The folowing vector spaces have the same dimension :
  \begin{itemize}
  \item $H^0(\SUtpar{2, K_X(D)}, \Ltpar{r})$
  \item $H^0(\mathcal M'(1,\vect\lambda^T), \mathcal L^{r,\vect\lambda^T})$
  \item $H^0(\mathcal M(r,\vect\lambda), \mathcal L^{1,\vect\lambda})$
  \item $H^0(\mathcal U_\Xp(r,\vect\lambda), \mathcal L^{2,\vect\lambda})$
  \end{itemize}
\end{prop}

\begin{proof}
  The isomorphism $Sp_2 \simeq SL_2$ gives the first equality, the second
  follows from the rank-level symmetry for Verlinde numbers of the
  symplectic group, and the last equality folows from the $SL_2/GL_n$
  rank-level duality. 
\end{proof}

\begin{thm}
  There is a commutative diagram of isomorphisms :
  \[ \xymatrix{
    H^0(\SUtpar{2, K_X(D)}, \Ltpar{r}) \ar[r]
    & H^0(\mathcal M'(1,\vect\lambda^T), \mathcal L^{r,\vect\lambda^T}) \\
    H^0(\mathcal U_\Xp(r,\vect\lambda), \mathcal L^{2,\vect\lambda})^\vee
    \ar[r] \ar[u]
    & H^0(\mathcal M(r,\vect\lambda), \mathcal L^{1,\vect\lambda})^\vee
    \ar[u] \\
  } \]
\end{thm}

\begin{proof}
  The bottom line comes from the inclusion $E \to E \oplus E^\vee$ between
  the moduli spaces, the top line is the identification of $SL_2$ with
  $Sp_2$, the left edge is the rank-level isomorphism. This implies the
  arrow from $H^0(\mathcal U_\Xp(r,\vect\lambda), \mathcal
  L^{2,\vect\lambda})^\vee$ to $H^0(\mathcal M'(1,\vect\lambda^T), \mathcal
  L^{r,\vect\lambda^T})$ is an isomorphism, so the two other arrows are
  isomorphisms.
\end{proof}


\end{document}